\documentstyle[draft]{article}
\title{Differential geometry of generalized almost quaterniohic structures, II}
\author{ V.F. Kirichenko, O.E. Arseneva}
\newcommand \implies{\Longrightarrow}
\newcommand \eq{\!\!&\!\!=\!\!&\!\!}
\newcommand \la{\langle}
\newcommand \ra{\rangle}
\newcommand \n[2]{\nabla_#1\,#2}

\newcommand \hq{{\cal H}AQ_\a}
\newcommand \kq{{\cal K}AQ_\a}
\newcommand \ah{A{\cal H}_\a}
\newcommand \k{{\cal K}_\a}

\newcommand \tr{\mbox{\,\rm tr\,}}
\newcommand \End{\mbox{\rm End\,}}
\newcommand \id{\mbox{\,\rm id}}
\newcommand \re{\mbox{\rm Re\,}}
\newcommand \im{\mbox{\rm Im\,}}
\newcommand \procl[1]{\medskip\par{\bf#1.}\it}
\newcommand \eprocl{\medskip\par\rm}
\newcommand \definition[1]{\medskip\par{\bf#1.}}
\newcommand \edefinition{\medskip\par}
\newcommand \example[1]{\medskip\par{\bf#1.}}
\newcommand \eexample{\medskip\par}
\newcommand \remark[1]{\medskip{\it#1.}}
\newcommand \eremark{\par}
\newcommand \demo[1]{\par{\it#1.}}
\newcommand \edemo{\medskip\par}
\newcommand \qed{$\quad\Box$}
\newcommand \qeds{\quad\Box}
\renewcommand \a{\alpha}
\newcommand \e{\varepsilon}
\renewcommand \b{\beta}
\renewcommand \c{\gamma}
\renewcommand \d{\delta}
\renewcommand \l{\lambda}
\renewcommand \L{\Lambda}
\newcommand \om{\omega}
\newcommand \z{\zeta}
\renewcommand \v{\vartheta}
\renewcommand \o{\overline}
\begin{document}
\maketitle
\setcounter{equation}{37}
\setcounter{section}{3}
\begin{abstract} The generalized conseption of self-dual and anti-self dual forms and
manifolds is built. The complete classification of important classes of
self-dual and anti-self-dual generalized Kaehler manifolds is obtained.
\end{abstract}
\tableofcontents
\newpage
The present paper is the immediate continuation of Part 1 (see page 00-00).
Herein we shall refer to the results of Paper 1. The main purpose of the present
paper is studying geometry originality of 4-dimensional (pseudo-) Riemannian
manifolds. This originality is contained in the existence of a self-dual
structure on such manifolds. Then we investigate generalization possibilities
of the structure on the manifolds of greater dimensions. Such investigation is
created on the base of the concept of generalized almost quaternionic structures
which was developed in the previous paper.
\par
Section 4 introduces and investigates the class of vertical type $AQ_\a$-struc\-tures
generalizing the class of quaternionic-Kaehler structures, and finds the criterion
of the structural bundle of the structure being Einsteinian that generalized
the well-known Atiyah-Hitchin-Singer criterion of 4-dimensional oriented
Riemannian manifolds being Einsteinian in terms of self-dual forms on the
manifolds. We also prove the theorem of generalized quaternionic-Kaehler manifolds
of dimension greater than 4 being Einsteinian that generalizes the well-known
result of M.Berger.
\par
Section 5 introduces and investigates the notion  of $t$-conformal-semiflat
$AQ_\a$-manifold generalizing the classical notion of self-dual and anti-self-dual
4-dimen\-sional Riemannian manifolds in the case $AQ_\a$-manifolds ofarbitrary
dimension. We introduce the notion of twistor curvature tensor of a vertical
type $AQ_\a$-manifold and prove that the twistor curvature tensor of
$t$-conformal-semiflat $AQ_\a$-manifold is an algebraic curvature tensor. We
also prove that anti-self-dual $AQ_\a$-manifolds as well as generalized
quaternionic-Kaehler manifolds are the manifolds of pointwise constant twistor
curvature. It is proved that generalized hyper-Kaehler manifolds are zero
twistor curvature manifolds, and they also are zero Ricci curvature manifolds.
The result generalizes the well-known result of M.Berger.
\par
In Section 6 we prove that a 4-dimensional generalized Kaehler manifold is
self-dual iff its (generalized) Bochner curvature tensor is equal to zero.
Moreover, we get a complete classification of 4-dimensional self-dual
non-exceptional manifolds of constant scalar curvature. It is proved that a
4-dimensional generalized Kaehler manifold is anti-self-dual if and only if
its scalar curvature is equal to zero. It is also proved that a 4-dimensional
compact regular spinor manifold carrying a classical type Kaehler structure is
anti-self-dual if and only if it is Ricci-flat. The above results essentially
generalize the well-known results of N. Hitchin, J.-P. Bourguignon, A. Derdzinski,
B.-Y. Chen and M. Itoh.
\section{Generalized Almost Quaternionic Manifolds of Vertical Type}
The study of geometry of 4-dimensional (pseudo-) Riemannian manifolds is of
great importance in the modern geometric investigations since such manifolds
are significant in mathematical physics. It was the features of the geometry
that gave rise to the wide use of the well-known Penrose Twistor Programme in
modern investigations in the theory of gravitation and Yang-Mills fields. These
features are determined to a great extent by the existence of a natural almost
quaternionic structure on 4-dimensional (pseudo-) Riemannian manifolds. Na\-turally
the que\-stion arises: what geometrical features of 4-dimensional manifolds admit
generali\-zation onto almost quaternionic manifolds of arbitrary dimension ?
For example, it is known that with every almost quaternionic manifold there is
naturally associated a twistor bundle, and a number of important results of
Penrose twistor geometry are naturally, though non-trivially, generalized on
other types of almost quaternionic manifolds, for example, on
quaternionic-Kaehler ones [1],[2].
\par
The present chapter distinguishes the class of generalized almost quater\-nionic
structures for which similar extrapolation looks most natural. They are the
so-called {\it generalizeg almost quaternionic structures of vertical type\/}
that, as we shall see below, generalize quaternionic-Kaehler structures. The
main result of the chapter shows that structural bundle of the manifold
carrying the $AQ_\a$-structure of vertical type is Einsteinian if and only if
its Kaehler module is invariant with respect to Riemann-Christoffel endomorphism.
This widely generalizes the classical Atiyah-Hitchin-Singer result giving a
criterion of 4-dimensional Riemannian manifolds being Einsteinian in terms of
self-dual forms [3]. It is also proved that generalized quaternionic-Kaehler
manifold of dimension greater than 4 is an Einsteinian manifold, that generalizes
the well-known Berger theorem of quater\-nionic-Kaehler manifolds being Einsteinian [1].
\subsection{Self-dual and anti-self-dual forms on $AQ_\a$-manifolds}
Let $M$ be an $AQ_\a$-manifold, $\sf Q$ be its structural bundle. Note that
algebra ${\bf H}_\a$ of generalized quaternions admits a canonical representation
(of the first kind) $\l:{\bf H}_\a\rightarrow \End{\bf H}_\a,\
\l(q)(X)=qX;\ X\in{\bf H}_\a$. Indeed, $\l(q_1\,q_2)(X)
=(q_1\,q_2)X=q_1(q_2\,X)=\l(q_1)\circ\l(q_2)(X)$, i.e. $\l(q_1\,q_2)
=\l(q_1)\circ\l(q_2);\quad q_1,q_2\in{\bf H}_\a$. Let $\{\,1,j_1,j_2,j_3\,\}$
be orthonormalized basis of the space ${\bf H}_\a,j_3=j_1j_2,q=a+bj_1+cj_2
+dj_3\in{\bf H}_\a$. Then in the basis
$$(\l(q)^i{}_j)=\left(\begin{array}{rrrr}
a&\a b&\a c&-d\\b&a&\a d&-\a c\\c&-\a d&a&\a b\\d&-c&b&a\end{array}\right).$$
\par
Thus, any $\a$-quaternion admits identification with endomorphism of 4-dimen\-sional
real space ${\bf H}_\a$. Lowering index of the endomorphism we get a tensor of
the type (2,0) on this linear space that will be skew-symmetric iff $a=0$, i.e.
$\a$-quaternion $q$  is purely imaginary. In its turn, fibre bundle $\sf Q$
has a natural metric $g=(\cdot,\cdot)$ generated by the metric of algebra
${\bf H}_\a$. Using it we lower index of the purely imaginary $\a$-quaternion
regarded in it canonical representation and get a 2-form $\om$ on the space of
fibre bundle $\sf Q$, that we call {\it a fundamental form of the\/} $\a${\it-quaternion}.
\par
Similarly, we can consider canonical representation of the second kind
$\mu:{\bf H}_\a\rightarrow \End{\bf H}_a;\ \mu(q)(X)=Xq;\ X\in{\bf H}_\a$.
In this case $\mu(q_1\,q_2)=\mu (q_2)\circ\mu(q_1);\quad q_1,q_2\in{\bf H}_\a$.
Here the mentioned $\a$-quaternion $q$ in orthonormalized basis $\{\,1,j_1,j_2,j_3\,\}$,
$j_3=j_1j_2$, will be represented by the matrix
$$(\mu(q)^i{}_j)=\left(\begin{array}{rrrr}
a&\a b&\a c&-d\\b&a&-\a d&\a c\\c&\a d&a&-\a b\\d&c&-b&a\end{array}\right)$$
Lowering index of endomorphism $\mu(q)$ we get the tensor of the type (2,0) on
the real space ${\bf H}_\a$ that will be skew-symmetric iff $a=0$. Thus, lowering
index of the purely imaginary $\a$-quaternion $q\in\{{\sf Q}\}$  regarded
in canonical representation of the second kind and using the fibre metric we
get 2-form on the space of fibre bundle $\sf Q$ that we call {\it a pseudo-fundamental
form of the} $\a${\it-quaternion.}
\par
Similarly to the above, an $AQ_\a$-manifold $M$ is assumed calibrated and thus,
${\sf X}({\sf Q})={\cal V}\oplus{\cal H}$, where $\cal V$ and $\cal H$ are,
respectively, vertical and horisontal distributions of metric connection in
fibre bundle $\sf Q$ induced by calibration. Recall that tensor $t\in{\cal T}^s_r({\sf Q})$
is called {\it vertical\/} if it vanishes if at least one of the arguments of
the tensor regarded as a polylinear function is horisontal. In other word, a
tensor $t$ is vertical if $\forall\,X_1,\ldots,X_r\in{\sf X}({\sf Q})$,
$\forall\,\om^1,\ldots,\om^s\in{\sf X}^*({\sf Q})\implies t(X_1,\ldots,X_r,
\om^1,\ldots,\om^s)=t(X^V_1,\ldots,X^V_r,\om^1_V,\ldots,\om^s_V)$, where
$X^V$, $\om_V$ are vertical constituents of vector $X$ and covector $\om$,
respectively. Vertical tensors define a subalgebra of tensor algebra ${\cal T}({\sf Q})$
that we call {\it vertical tensor algebra\/} and denote it by ${\cal T}_V({\sf Q})$.
Similar remarks refer to Grassnanian algebra $\L({\sf Q})\supset\L_V({\sf Q})$.
In particular, pointwise localization of module $(\L^r)_V({\sf Q})$ has dimension
$C^r_4=\frac{4!}{r!(4-r)!};\quad 0\le r\le4$, and zero dimension, when $r>4$.
Thus, we have naturally defined the Grassmanian algebra $\L_V({\sf Q})
=\oplus_{r=0}^4(\L^r)_V({\sf Q})$ that we call {\it a vertical Grassmanian
algebra.} Note that module $(\L^4)_V$ is one-dimensional; as its basis we
naturally take the form $\eta$ defined by the equality
$$\eta_p=\sqrt{\det(g)}\,\om^0\land\om^1\land\om^2\land\om^3,$$
where $(g)$ is Gramme matrix of metric $g_p$, $\{\,p,\om^0,\om^1,\om^2,\om^3\,\}$
is a coframe dual to frame $\{\,p,e_0,e_1,e_2,e_3\,\}\in B{\sf Q}$ of the space
${\cal V}_p={\sf Q}_p;\quad p\in{\sf Q}$. We check the independence of $\eta_p$
on the choice of basis in a standard way. We call the form $\eta$ {\it a vertical
volume form}. Evidently, it generates the orientation of the bundle ${\cal V}=\{{\sf Q}\}$.
\definition{Definition 18} We call {\it (vertical) Hodge operator} the endomorphism
$*:\L_V({\sf Q})\rightarrow\L_V({\sf Q})$ being defined by the equality
$$\om\land(*\v)=(\om,\v)\eta;\quad\om,\v\in(\L^r)_V({\sf Q}),\quad r=0,1,2,3,4;$$
where $(\cdot,\cdot)$ is scalar product in $\L_V({\sf Q})$ induced by metric in $\sf Q$.
\edefinition
In particular, $*:(\L^2)_V({\sf Q})\rightarrow(\L^2)_V({\sf Q})$, and $(\om,\v)
=\frac12\om_{\b\c}\,\v^{\b\c}$, then $(*\om)_{\b\c}=\om_{\d\e}$, where
$\left(\begin{array}{cccc}0&1&2&3\\\b&\c&\d&\e\end{array}\right)$ is even
permutation. Hence, $(*)^2=\id$ on $(\L^2)_V({\sf Q})$ and then,
$*|_{(\L^2)_V({\sf Q})}$ has eigenvalues $\pm1$.
\definition{Definition 19} A vertical 2-form on $\sf Q$ is called
{\it self-dual\/} (resp., {\it anti-self-dual\/}) if it is an eigenvector of
Hodge operator with eigenvalue 1 (resp., -1). If it is also parallel along the
fibres of structural bundle (i.e. projectable) it is called {\it a form on M.}
The module of self-dual (resp., anti-self-dual) forms on $M$ will be denote by
$\L^+(M)$ (resp., $\L^-(M))$.
\edefinition
\procl{Theorem 23} A projectable vertical 2-form is self-dual on an
$AQ_\a$-manifold $M$ iff it is a fundamental form of a purely imaginary
$\a$-quaternion. Here, the orientation of fibre bundle $\sf Q$ generated by
exterior square of the form corresponds to the canonical orientation of the
bundle iff the $\a$-quaternion has a real norm.
\eprocl
\demo{Proof} Let $\om$ be a projectable vertical 2-form. By definition,
$\om\in\L^+(M)$ $\iff*\om=\om$. Denote by $\{\,\eta_{\b\c\d\e}\,\}$ components of
tensor $\eta$ in a positively oriented frame. Then, if $\{\,\om_{\b\c}\,\}$ are
components of the 2-form $\om$, then the condition of its self-duality (see (7))
can be written in the form
$$\om_{\b\c}=\frac12\eta_{\b\c\d\e}\,g^{\d\z}\,g^{\e\v}\,\om_{\z\v}.$$
On the space of bundle $B{\sf Q}$, where $(g_{\b\c})=diag(1,-\a,-\a,1)$ these
correlations will assume the form $\om_{\hat\b\hat\c}=\e(\b,\c)\,\om_{\b\c}$,
where $\e(\b,\c)=g_{\b\b}\,g_{\c\c};\quad\b,\c=0,1,2,3$; $(\hat\b,\hat\c)$
is a complementing the pair $(\b,\c)$ up to the even permutation of indices
(0,1,2,3). It means that
\begin{equation}
\om\in\L^+(M)\iff(\om_{\b\c})=\left(\begin{array}{rrrr}
0&\a x&\a y&-z\\-\a x&0&-z&y\\-a y&z&0&-x\\z&-y&x&0\end{array}\right)
\end{equation}
Raising the index of the form we get the endomorphism of the fibre bundle
$\sf Q$ that, by the above, can be identified with the purely imaginary
$\a$-quaternion $q=xJ_1+yJ_2+zJ_3$. Evidently, the inverse is also true: by (38)
the fundamental form of purely imaginary $\a$-quaternion is self-dual. Note
that, by definition of self-dual form, $\om\land\om=\om\land(*\om)=(\om,\om)\eta$,
and hence, orientation of the fibre bundle $\sf Q$ generated by the form
$\om\land\om$ coincides with the canonical orientation of this fibre bundle iff
$\|\om\|^2>0$, i.e. $|q|^2>0$.\qed
\edemo
Similarly, we prove
\procl{Theorem 24} A projectable vertical 2-form is anti-self-dual on a
$AQ_\a$-manifold iff it is a pseudofundamental form of a purely imaginary
$\a$-quaternion. Here, orientation of the fibre bundle generated by exterior
square of the form is inverse to the natural orientation of the fibre bundle
iff the $\a$-quaternion has a real norm.\qed
\eprocl
\subsection{Vertical tensors on $AQ_\a$-manifold}
Let $M$ be an $4n$-dimensional almost $\a$-quaternionic manifold, $\sf Q$ be
its structural bundle. In a standard way [4] it is proved that giving of an
$AQ_\a$-structure on $M$ is equivalent to giving of $\cal G$-structure on the
manifold with structural group ${\cal G}=GL(n,{\bf H}_\a)\cdot Sp_\a(1)$, where
$Sp_\a(m)$ is a symplectic group of order $m$ over the ring ${\bf H}_\a$. The
elements of the space of this $\cal G$-structure called {\it adapted frames},
or {\it A-frames}, are constructing in the following way. Let $p\in M$,
$r=(J_0=\id,J_1,J_2,J_3)\in B{\sf Q}$ be a positively oriented orthonormalized
basis of fibre ${\sf Q}_p$, $(e_1,\ldots,e_{n})$ be the basis of the space
$T_p(M)$ regarded as ${\bf H}_\a$-module: if $q=a+bj_1+cj_2+dj_3\in{\bf H}_\a,
\quad X\in T_p(M)$, then $qX=aX+bJ_1X+cJ_2X+dJ_3X$. Assume the greek indices to
be in range from 0 to 3, the latin ones - in range from 1 to $n$. Denote
$e_{\b a}=J_\b(e_a)$. Then $(p,e_{\b a};\quad\b=0,\ldots,3;\quad a=1,\ldots,n)$
is a frame of the space $T_p(M)$ called {\it adapted.} Evidently, in this frame
\begin{eqnarray}
(J_1)=\left(\begin{array}{rrrr}
0&\a I&0&0\\I&0&0&0\\0&0&0&\a I\\0&0&I&0\end{array}\right);\qquad
(J_2)=\left(\begin{array}{rrrr}
0&0&\a I&0\\0&0&0&-\a I\\I&0&0&0\\0&-I&0&0\end{array}\right);\nonumber\\
(J_3)=\left(\begin{array}{rrrr}
0&0&0&-I\\0&0&\a I&0\\0&-\a I&0&0\\I&0&0&0\end{array}\right).
\hphantom{---------}\end{eqnarray}
and then $(J_k)=(j_k)\otimes I;\quad k=1,2,3;\quad I$ is a unit matrix of order
$n$. Evidently, giving tensor $t$ of the type (1,1) on $M$ is defined by giving
of a function set $\{\,t^{\b i}_{\c j}\,\}$ on the space of $\cal G$-structure
that are components of the tensor in the corresponding A-frame. It follow from
(39) that such tensor is an $\a$-quaternion $q=a\id+bJ_1+cJ_2+dJ_3$ iff
$t^{\b i}_{\c j}=q^\b_\c\,\d^i_j$, where
$$(q^\b_\c)=\left(\begin{array}{rrrr}
a&\a b&\a c&-d\\b&a&\a d&-\a c\\c&-\a d&a&\a b\\d&-c&b&a\end{array}\right)$$
is the matrix of $\a$-quaternion $q$ in basis $\{\id,J_1,J_2,J_3\,\}$
regarded as endomor\-phism of the structural bundle.
\par
Consider the natural representation of endomorphisms algebra of module ${\sf X}(M)$
of $AQ_\a$-manifold $M$ into endomorphisms algebra of module ${\cal T}^1_1(M)$,
the represen\-tation being generated by the left shifts. Namely, if $f\in{\cal T}^1_1(M)$
we juxtapose it to endomorphism $\hat f$ acting by formula $\hat f(g)=f\circ g;
\quad g\in{\cal T}^1_1(M)$. Evidently, the mapping $f\rightarrow\hat f$ is
representation in view of endomorphisms algebra being associative. The question
naturally arises: when does element $f\in{\cal T}^1_1(M)$ in this representation
preserves the structural bundle and, thus, induces the structural bundle
endomorphism ? The question is answered by
\procl{Theorem 25} Let $M$ be an $AQ_\a$-manifold. Endomorphism $t\in{\cal T}^1_1(M)$
induces structural bundle endomorphism iff its components on the space of
$\cal G$-structure have the form
\begin{equation}
t^{\c c}_{\b b}=t^\c_\b\,\d^c_b.
\end{equation}\eprocl
\demo{Proof} Let equalities (40) hold and let $q\in\{\,{\sf Q}\,\}$. Assume
$\hat t(q)=t\circ q$. Fix frame $(p,J_\b;\quad\beta=0,\ldots,3)$ in ${\sf Q}_p$
and its corresponding $A$-frame $(p,e_{\b b})\in{\cal G}$. We have
$\hat t(J_\b)(e_b)=t\circ J_\b(e_b)=t(J_\b(e_b))=t(e_{\b b})=t^{\c c}_{\b b}\,e_{\c c}
=t^\c_\b\,\d^c_b\,e_{\c c}=t^{\c}_{\b}\,e_{\c b}=t^\c_\b\,J_\c(e_b);\quad
b=1,\ldots,n$, hence $\hat t(J_\b)=t^\c_\b\,J_\c$, i.e. endomorphism
$\hat t: f\rightarrow t\circ f$ of module ${\cal T}^1_1(M)$ preserves the structural
bundle, and thus, induces its endomorphism. Inversely, let $\hat t:
f\rightarrow t\circ f$ preserves the structural bundle. Then $\hat t(J_\b)=t^\c_\b\,J_\c$
and $t(e_{\b b})=t(J_\b(e_b))=t\circ J_\b(e_b)=\hat t(J_\b)(e_b)=t^\c_\b\,J_\c(e_b)
=t^\c_\b\,\d^c_b\,J_\c(e_c)=t^\c_\b\,\d^c_b\,e_{\c c}$. Thus, $t^{\c c}_{\b b}=t^\c_\b\,\d^c_b$.\qed
\edemo
Evidently, giving structural bundle endomorphism is equivalent to giving vertical
tensor $t_V\in({\cal T}^1_1)_V({\sf Q})$ and, thus, Theorem 25 can be formulated as follow:
\procl{Theorem 26} Let $M$ be an $AQ_\a$-manifold. Endomorphism
$t\in{\cal T}^1_1(M)$ in natural representation generates vertical tensor
$t_V\in({\cal T}^1_1)_V({\sf Q})$ iff its compo\-nents on the space of $\cal G$-structure
have the form $t^{\c c}_{\b b}=t^\c_\b\,\d^c_b$.\qed
\eprocl
\definition{Definition 20} Endomorphism $t$ of the module ${\sf X}(M)$ of a
$AQ_\a$-manifold $M$ preserving the structural bundle in the natural representation
and, thus, generating vertical tensor on $\sf Q$ is called {\it vertical endomorphism.}
\edefinition
\example{Example} As we have seen, any $\a$-quaternion $q\in\{\,{\sf Q}\,\}$
satisfies the condition of Theorem 26 and, thus, it is a vertical endomorphism.
The validity of (40) in this case can also be shown in another way. Let $q\in{\sf Q}$.
Then $q(e_{\b b})=q(J_\b(e_b))=q^\c\,J_\c(J_\b(e_b))=q^\c(J_\c\circ J_\b)e_b
=q^\c\,C^\e_{\c\b}\,J_\e(e_b)=q^\c\,C^\e_{\c\b}\,\d^c_b\,J_\e(e_c)=q^\e_\b\,\d^c_b\,e_{\e c}$
i.e. $q^{\e c}_{\b b}=q^\e_\b\,\d^c_b$, where $q^\e_\b= C^\e_{\c\b}\,q^\c,
\{\,C^\e_{\c\b}\,\}$ are structural tensor components of $\a$-quaternion algebra.
\eexample
Similar terminology will be preserved for other tensor types on $M$, received
from vertical endomorphisms by means of classical operations of tensor algebra.
In particular, if $M$ is a pseudo-Riemannian manifold, a 2-form on $M$ is called
{\it vertical\/} if it is corresponds to a vertical endomorphism in raising
the index. Evidently, the set of all vertical 2-forms on $M$ forms a submodule
$(\L^2)_V(M)\subset\L^2(M)$ that we call {\it a module of vertical 2-forms on M.}
\definition{Definition 21} An $AQ_\a$-structure on a pseudo-Riemannian manifold
$(M,G=\la\cdot,\cdot\ra)$ is called {\it generalized  quaternionic-Hermitian},
or $\hq${\it-structure} if
$$\forall J\in\{\,{\sf T}\,\}\implies\la JX,Y\ra+\la X,JY\ra=0;\qquad X,Y\in{\sf X}(M).$$
\edefinition
In this case for any purely imaginary $\a$-quaternion $q$ on $M$ tensor
$\Omega(X,Y)=\la X,qY\ra$ is a vertical form on $M$ called {\it Kaehler form\/} of
$\a$-quaternion $q$. Kaehler 2-forms on $M$ form a submodule
${\cal K}(M)\subset(\L^2)_V(M)$ that we call {\it Kaehler module\/} of $\hq$-structure.
By Theorem 23, the image of Kaehler module in natural representation (in
combination with lowering the index by means of fibre metric) is a module of
self-dual forms on $M$.
\par
Giving an $\hq$-structure on manifold $M^{4n}$ is equivalent to giving
$\cal G$-structure on $M$ with structural group ${\cal G}=Sp_\a(n)\cdot Sp_\a(1)$
whose total space elements are $A$-frames, and their vectors $\{\,e_1,\ldots,e_n\,\}$
form an orthonormalized system in quaternionic-Hermitian metric $\la\la X,Y\ra\ra
=\la X,Y\ra+i\la X,I^3Y\ra+j\la X,J^3Y\ra+k\la X,K^3Y\ra$, where
$\{\id,I,J,K\,\}$ is positively oriented orthonor\-malized basis of fibre
${\sf Q}_p$ an arbitrary point $p\in M$.
\procl{Lemma 3} On the space of $\cal G$-structure components of metric
tensor of an $\hq$-manifold have the form: $G_{\b b\c c}=g_{\b\c}\,G_{bc}$, where
$\{\,g_{\b\c}\,\}$ are components of fibre metric of the structural bundle,
$G_{bc}=\la e_b,e_c\ra$.
\eprocl
\demo{Proof} We have:
\begin{eqnarray*}
G_{\b b\c c}\eq\la J_\b(e_b),J_\c(e_c)\ra=-\la J_\c\circ J_\b(e_b),e_c\ra=\\
\eq-C^\d_{\c\b}\la J_\d(e_b),e_c\ra=-C^0_{\c\b}\la e_b,e_c\ra,
\end{eqnarray*}
where $\{\,C^\d_{\c\b}\,\}$ are structural constants of $\a$-quaternion
algebra. Note that since for any pair $\{\,q_1,q_2\,\}$ of purely imaginary
$\a$-quaternions $q_1q_2+q_2\,q_1=-2\la q_1,q_2\ra,\ q_1q_2-q_2\,q_1=[q_1,q_2]$,
then $q_1\,q_2=-\la q_1,q_2\ra+\frac12[q_1,q_2]$. Thus, $\re(q_1q_2)=-\la q_1,q_2\ra$,
$\im(q_1q_2)=\frac12[q_1,q_2]$, where $\re q$ and $\im q$ are real and purely
imaginary parts of $\a$-quaternion $q$, respectively. Hence, $C^0_{\c\b}
=\re(J_\c J_\b)=-\la J_\c,J_\b\ra=-g_{\c\b}=-g_{\b\c}$. Thus, $G_{\b b\c c}
=g_{\b\c}\la e_b,e_c\ra=g_{\b\c}\,G_{bc}$.\qed
\edemo
Let $\Omega$ be a vertical 2-form on an $\hq$-manifold $M$, $\{\,\Omega_{\b b\c c}\,\}$
be its components on the space of $\cal G$-structure. By Theorem 26,
$G^{\e h\b b}\,\Omega_{\e h\c c}=t^\b_\c\,\d^b_c$, and thus,
$G_{\b b\d d}\,G^{\e h\b b}\,\Omega_{\e h\c c}=t^{\b}_{\c}\,G_{\b c\d d}$. But
$G_{\b b\d d}\,G^{\e h\b b}=\d^\e_\d\,\d^h_d$, hence, $\Omega_{\d d\c c}
=G_{\d d\b c}\,\d^\b_\c$. But by Lemma 3,
\begin{equation}
\Omega_{\d d\c c}=\om_{\d\c}\,G_{dc},
\end{equation}
where $\om_{\d\c}=g_{\d\b}\,t^{\b}_{\c}$ are components of tensor $\om$ on
$\sf Q$, being the image of form $\Omega$  in natural representation. From
(41) it follows that $\om_{\d\c}=\frac1nG^{dc}\,\Omega_{\d d\c c}$, in particular,
$\om$ is skew-symmetric tensor being, evidently, a projectable vertical 2-form
on the space $\sf Q$. Evidently, the inverse is also true: let $\om$ be
projectable vertical 2-form on $\sf Q$, then 2-form $\Omega$ on $M$ defined
by (41) is a vertical 2-form on $M$. Thus, we have proved
\procl{Theorem 27} The module of vertical 2-forms on an $\hq$-manifold $M$
in natural representation coincides with module $(\L_{\pi V})^2({\sf Q})$ of
projectable vertical 2-forms on the space $\sf Q$.\qed
\eprocl
Now we shall identify the moduls sometime.
\par
As a corollary we get the following result:
\procl{Theorem 28} The module of vertical 2-forms on an $\hq$-manifold $M$
is decom\-posed into orthogonal direct sum of Kaehler module and the submodule
coinciding in natural representation with the submodule of anti-self-dual forms
on $M$ and serving as a Kaehler module of a uniquely defined $\hq$-structure on M.
\eprocl
\demo{Proof} Note that module $(\L_{\pi V})^2({\sf Q})$ of projectable vertical
2-forms on $\sf Q$ is decom\-posed into orthogonal direct sum of submodules of
self-dual and anti-self-dual forms on $M$. Indeed, by definition,
$(\L_{\pi V})^2({\sf Q})=\L^+(M)\oplus\L^-(M)$.If $\om\in\L^+(M),\v\in\L^-(M)$,
then $\la\om,\v\ra\eta=\om\land(*\v)=-\om\land\v=-\v\land\om=-\v\land(*\om)
=-\la\v,\om\ra\eta=-\la\om,\v\ra\eta$, hence, $\la\om,\v\ra=0$, and thus,
$\L^+(M)\perp\L^-(M)$. Now let $\Omega,\Theta\in\L_V(M)$, $\om,\v$ are their
images in natural representation. From (41) it follows that the forms $\Omega$
and $\Theta$ are orthogonal in metric $G$ iff the forms $\om$ and $\v$ are
orthogonal in fibre metric. In particular, submodule ${\cal K}(M)^\perp$ in
natural representation coincides with submodule of anti-self-dual forms on $M$.
Further, as in proof of Theorem 4, we see that in natural identification of
elements of module ${\cal K}(M)^\perp$ with endomorphisms of module ${\sf X}(M)$
fibre bundle ${\cal S}c(M)\oplus{\cal K}(M)^\perp$ defines an $\hq$-structure on $M$.\qed
\edemo
\definition{Definition 22} $\hq$-structures ${\sf Q}={\cal S}c(M)\oplus{\cal K}(M)$
and ${\sf Q}^*={\cal S}c(M)\oplus{\cal K}(M)^\perp$ on manifold $M$ will be called {\it conjugate.}
\edefinition
\subsection{Generalized quaternionic-Kaehler structures}
\definition{Definition 23} An $\hq$-structure is called {\it a generalized
quaternionic-Kaehler}, or $\kq$-structure if its structural bundle is invariant
with respect to parallel translations in Riemannian connection.
\edefinition
It means that Riemannian connection $\nabla$ is an $AQ_\a$-connection, and we
fix it as calibration.
\par
Evidently, any quaternionic-Kaehler structure (see example 3 in section 2.2)
is a $\kq$-structure.
\procl{Theorem 29} A $kq$-manifold $M$ of dimension $4n>4$ is an Einsteinian
manifold.\eprocl
\demo{Proof} Let ${\sf U}=\{\,U_a\,\}_{a\in A}$ be the local triviality covering
of $\cal G$-structure for manifold $M$, $U\in{\sf U}$. Then narrowing $\sf Q$
on $U$ is a $\pi AQ_\a$-structure, and thus, there is a pair twistors $\{\,I,J\,\}$
on $U$, such that $\{\,{\sf Q}|U\,\}$ is generated by their algebraic shell,
i.e. $\{\,{\sf Q}|U\,\}={\cal L}(\id,J_1,J_2,J_3\,)$, where $J_1=I,J_2=J,J_3=K= I\circ J$.
Besides, there exists the system $\{\,e_1,\ldots,e_n\,\}$ of vector fields on
$U$ such that the system $\{\,e_{\b b}\,\},\quad\b=0,\ldots,3;\quad b=1,\ldots,n$,
where $e_{\b b}=J_\b(e_b)$, defines the space section of $\cal G$-structure
over $U$. In this case if $X,Y\in{\sf X}(U)$ then $\n Y(J_\b)=\om(Y)^\c_\b\,J_\c$,
where $\om=\{\,\om^\c_\b\,\}$ is a form on $U$ with values in Lie algebra
${\sf so}(2-\a,1+\a;{\bf R})$ of the structural group of bundle $B{\sf Q}$,
at every point of $U$. Further, let $R(X,Y)=\n X\n Y-\n Y\n X-\nabla_{[X,Y]}$
is a Riemann-Christoffel tensor. Then
\begin{eqnarray*}
R(X,Y)(J_\b)\eq\n X\n Y(J_\b)-\n Y\n X(J_\b)-\nabla_{[X,Y]}(J_\b)\\
\eq\n X(\om(Y)^\c_\b\,J_\c)-\n Y(\om(X)^\c_\b\,J_\c)-\om([X,Y])^\c_\b\,J_\c\\
\eq X\om(Y)^\c_\b\,J_\c+\om(Y)^\c_\b\,\n XJ_\c-Y\om(X)^\c_\b\,J_\c-\\
&{}&-\om(X)^\c_\b\,\n YJ_\c-\om([X,Y])^\c_\b\,J_\c\\
\eq(X\om(Y)^\c_\b-Y\om(X)^\c_\b-\om([X,Y])^\c_\b)J_\c+\\
&{}&+\om(Y)^\c_\b\,\om(X)^\d_\c\,J_\d-\om(X)^\c_\b\,\om(Y)^\d_\c\,J_\d\\
\eq2d\om(X,Y)^\c_\b\,J_\c+(\om(X)^\d_\c\,\om(Y)^\c_\b-\om(Y)^\d_\c\om(X)^\c_\b)J_\d\\
\eq2(d\om(X,Y)^\c_\b+\frac12[\om,\om](X,Y)^\c_\b)J_\c.\end{eqnarray*}
Thus,
\begin{equation}
R(X,Y)J_\b=2D\om(X,Y)^\c_\b\,J_\c;\qquad\b,\c=1,2,3,
\end{equation}
where $D\om=d\om+\frac12[\om,\om]$ is a 2-form on $U$ with values in Lie algebra
${\sf so}(2-\a,$ $1+\a;{\bf R})$. Here $R(X,Y)$ is regarded as differentiation
of tensor algebra on manifold $M$ generated by endomorphism $R(X,Y)$ of module
${\sf X}(M)$, and thus,
\begin{equation}
R(X,Y)J_\b=[R(X,Y),J_\b].
\end{equation}
Let the form $D\om$ be given by the function matrix
$$(D\om^\c_\b)=\left(\begin{array}{rrr}
0&-c&b\\&0&-a\\\a b&-\a a&0\end{array}\right).$$
Then by (43) identities (42) can be written in the form
\begin{eqnarray}
1)\ [R(X,Y),I]\eq\hphantom{------}c(X,Y)J+\a b(X,Y)K;\nonumber\\
2)\ [R(X,Y),J]\eq-c(X,Y)I\hphantom{-----}-\a a(X,Y)K;\\
3)\ [R(X,Y),K]\eq b(X,Y)I-\ a(X,Y)J.\hphantom{------}\nonumber
\end{eqnarray}
Apply $(44_3)$ to vector $Z\in{\sf X}(U)$ and find scalar product of the
result and vector $JZ$:
\begin{eqnarray}
\la R(X,Y)KZ,JZ\ra-\la K\circ R(X,Y)Z,JZ\ra=\nonumber\\
=b(X,Y)\la IZ,JZ\ra-a(X,Y)\la JZ,JZ\ra.
\end{eqnarray}
Since $\la IZ,JZ\ra=-\la Z,KZ\ra=0$; $\la JZ,JZ\ra=-\a\la Z,Z\ra$;
$\la R(X,Y)KZ,JZ\ra=-\la R(X,Y)JZ,KZ\ra$; $-\la K\circ R(X,Y)Z,JZ\ra
=\la J\circ I\circ R(X,Y)Z,JZ\ra=-\a\la I\circ R(X,Y)Z,Z\ra=\a\la R(X,Y)Z,IZ\ra$,
equality (45) can be rewritten in the form
\begin{equation}
a(X,Y)\|Z\|^2=\la R(X,Y)Z,IZ\ra-\a\la R(X,Y)JZ,KZ\ra.
\end{equation}
Numerate elements $\{\,e_{\b b}\,\}$ of the local section of $\cal G$-structure
space by one index $\{\,e_i\,\},\quad i=1,\ldots,4n$. Assume $Z=\|e_i\|e_i$
in (33) and summarize by $i$ from 1 to $4n$. Remark that
$-\a\sum_{i=1}^{4n}\la R(X,Y)Je_i,Ke_i\ra\|e_i\|^2
=\sum_{i=1}^{4n}\la R(X,Y)Je_i,I\circ Je_i\ra\|Je_i\|^2
=\sum_{i=1}^{4n}\la R(X,Y)e_i,Ie_i\ra\|e_i\|^2$, and $\|e_i\|^4=1$. Then we get:
\begin{equation}
4na(X,Y)=2\sum_{i=1}^{4n}\la R(X,Y)e_i,Ie_i\ra\|e_i\|^2.
\end{equation}
Apply Ricci identity to the right part of the equality:
$$2na(X,Y)=-\sum_{i=1}^{4n}\{\la R(X,e_i)Ie_i,Y\ra+\la R(X,Ie_i)Y,e_i\ra\}
\|e_i\|^2,$$
or
$$2na(X,Y)=\sum_{i=1}^{4n}\{\la R(X,e_i)Y,Ie_i\ra-\la R(X,Ie_i)Y,e_i\ra\}
\|e_i\|^2.$$
Note that
\begin{eqnarray*}
\sum_{i=1}^{4n}\la R(X,Ie_i)Y,e_i\ra\|e_i\|^2=\sum_{i=1}^{4n}\la R(X,I^2e_i)Y,Ie_i\ra\|Ie_i\|^2\\
=-\sum_{i=1}^{4n}\la R(X,e_i)Y,Ie_i\ra\|e_i\|^2,\end{eqnarray*}
and then,
$$na(X,Y)=\sum_{i=1}^{4n}\la R(X,e_i)Y,Ie_i\ra\|e_i\|^2
=-\sum_{i=1}^{4n}\la I\circ R(X,e_i)Y,e_i\ra\|e_i\|^2.$$
But in view of $(44_1)$,
\begin{eqnarray*}
-\sum_{i=1}^{4n}\la R(X,e_i)Y,e_i\ra\|e_i\|^2=-\sum_{i=1}^{4n}\{\,\la R(X,e_i)IY,e_i\ra
-c(X,e_i)\la JY,e_i\ra-\\-\a b(X,e_i)\la KY,e_i\ra\,\}\|e_i\|^2=-Ric(X,IY)+c(X,JY)+\a b(X,KY),
\end{eqnarray*}
hence,
\begin{equation}
na(X,Y)=-Ric(X,IY)+c(X,JY)+\a b(X,KY).
\end{equation}
Where $Ric$ is Ricci tensor. Similarly,
\begin{eqnarray}
nb(X,Y)=-Ric(X,JY)-\a a(X,KY)-c(X,IY),\\
nc(X,Y)=\a\,Ric(X,KY)-\a b(X,IY)+\a a(X,JY).
\end{eqnarray}
Substituting $IY,JY$ and $KY$ for $Y$ in (48),(49) and (50), respectively, and
in view of $I^2=J^2=\a\id$, $K^2=-\id$, we get:
\begin{eqnarray}
na(X,JY)+b(X,JY)+c(X,KY)=-\a\,Ric(X,Y);\nonumber\\
a(X,IY)+nb(X,JY)+c(X,KY)=-\a\,Ric(X,Y);\\
a(X,IY)+b(X,JY)+nc(X,KY)=-\a\,Ric(X,Y).\nonumber
\end{eqnarray}
Substracting elementwise equations of the system, we get:
$$a(X,IY)=b(X,JY)=c(X,KY),$$
and in view of (51),
\begin{equation}
a(X,IY)=b(X,JY)=c(X,KY)=-\frac\a{n+2}\,Ric(X,Y).
\end{equation}
Now we can rewrite (46) in the form
$$\la R(X,IY)Z,IZ\ra-\a\la R(X,IY)JZ,KZ\ra=-\frac\a{n+2}\,Ric(X,Y)\|Z\|^2.$$
\noindent In particular, if $X=Y$, the equality has the form
\begin{equation}
\la R(X,IX)Z,IZ\ra-\a\la R(X,IX)JZ,KZ\ra=-\frac\a{n+2}\,Ric(X,X)\|Z\|^2.
\end{equation}
Substitute $JX$ for $X$ in (53):
\begin{equation}
\la R(JX,KX)Z,IZ\ra-\a\la R(JX,KX)JZ,KZ\ra=\frac\a{n+2}\,Ric(JX,JX)\|Z^2.
\end{equation}
Note that
$$Ric(JX,JX)=-\a(n+2)b(JX,J^2X)=-(n+2)b(JX,X)=-\a\,Ric(X,X).$$
In view of above multiply both parts of (54) by $(-\a)$ and add elementwise to (53):
\begin{eqnarray*}
\la R(X,IX)Z,IZ\ra-\a\la R(X,IX)JZ,KZ\ra-\a\la R(JX,KX)Z,IZ\ra+\\
+\la R(JX,KX)JZ,KZ\ra=-\frac{2\a}{n+2}\,Ric(X,X)\|Z\|^2.
\end{eqnarray*}
Note that the left part of the equality is symmetric with respect to $X$ and
$Z$. Thus, $Ric(X,X)\|Z\|^2=Ric(Z,Z)\|X\|^2$. Polarize this equality by argument $X$,
we get: $Ric(X,Y)\|Z\|^2=Ric(Z,Z)\la X,Y\ra$. Assume here $Z=\|e_i\|e_i$ and
summarize by $i$ from 1 to $4n$: $4n\,Ric(X,Y)=s\la X,Y\ra$, where
$s$ is scalar curvature of the manifold. Thus, $$Ric(X,Y)=\frac s{4n}\la X,Y\ra,$$
i.e. $M$ is Einsteinian manifold.\qed
\edemo
The above result generalizes the well-known Berger theorem of quaternionic-Kaehler
manifolds of dimension greater than 4 being Einsteinian [1].
\subsection{$\hq$-structures of vertical type}
Let $M$ be an $\hq$-manifold, $R$ be a Riemann-Christoffel tensor of metric
$G$. The classical features of symmetry of this tensor mean that it can be
regarded as self-conjugate endomorphism of a 2-form module on $M$. It will be
called a {\it Riemann-Cristoffel endomorphism}.
\definition{Definition 24} An $\hq$-structure on manifold $M$ is called
{\it the structure of vertical type} if a Riemann-Christoffel endomorphism
generated by the manifold metric preserves the module of vertical 2-forms on $M$.
\edefinition
For example, by Theorem 26 a canonical $\hq$-structure of a 4-dimensional oriented
pseudo-Riemannian manifold has a vertical type. Moreover, there holds the following
\procl{Theorem 30} Any $\kq$-manifold is a $\hq$-manifold of vertical type.
\eprocl
\demo{Proof} By definition, the structural bundle of $\kq$-manifold $M$, $\dim M=4n$,
is invariant with respect to parallel translations in Riemannian connection
$\nabla$ of manifold. The module of vertical endomorphisms as well as the one
of vertical 2-forms on $M$ have, naturally, the same property. (Indeed, if
$t\in{\cal T}^1_1(M)$ is a vertical endomorphism, $X\in{\sf X}(M)$, then
$(\n X\hat t)q=\n X(t\circ q)-t\circ \n X(q)\in\{\,{\sf Q}\,\};\quad q\in\{\,{\sf Q}\,\})$.
Moreover, in view of the parallelity
of the metric tensor in Riemannian connection the $\hq$-structure conjugate to
a $\kq$-structure is a $\kq$-structure itself. Now let $n>1$, $I$ be a twistor
at an arbitrary point $p\in M$, $\Omega$ be its Kaehler form. Then in view of
(47) and (52) we have
\begin{eqnarray*}
R(\Omega)(X,Y)\eq\frac12 R(X,Y)_{ij}\Omega^{ij}\\
\eq\frac12\sum_{i,j=1}^{4n}\la R(X,Y)e_i,e_\ra\Omega(e_i,e_j)
=\frac12\sum_{i,j=1}^{4n}\la R(X,Y)e_i,e_j\ra\la e_i,Ie_j\ra\\
\eq\frac12\sum_{i,j=1}^{4n}\la R(X,Y)e_i,Ie_j\ra\la e_i,I^2e_j\ra
=\frac\a 2\sum_{i=1}^{4n}\la R(X,Y)e_i,Ie_i\ra\|e_i\|^2\\
\eq\a na(X,Y)=na(X,I^2Y)=-\frac{\a n}{n+2}Ric(X,IY),
\end{eqnarray*}
where $Ric$ is Ricci tensor of manifold $M,\quad X,Y\in{\sf X}(M)$. Here
$\{\,e_1,\dots,e_{4n}\,\}$ is orthonormalized bazis adapted to twistor $J$.
But by Theorem 29, a $\kq$-manifold of dimension $4n>4$ is Einsteinian and thus,
$$R(\Omega)(X,Y)=-\frac{\a n\e}{n+2}\la X,IY\ra=-\frac{\a n\e}{n+2}\Omega(X,Y),$$
i.e.
\begin{equation}
{R(\Omega)=-\frac{\a n\e}{n+2}\Omega,}
\end{equation}
where $\e$ is Einstein constant. In particular, endomorphism $R$ transfers the
Kaehler module of a $\kq$-structure into itself. But in view of the same
considerations this endomorphism transfers into itself the Kaehler module of
the conjugate structure, too. And by Theoreme 28 the endomorphism preserves
the module of vertical 2-forms on $M$.\qed
\edemo
\procl{Corollary} Any quaternionic-Kaehler manifold is an $\hq$-manifold of
vertical type.\qed
\eprocl
Let $M$  be an $\hq$-manifold of vertical type. Then the Riemannian-Christoffel
endomorphism $R$ induces endomorphism $r$ of module $(\L^2)_V(M)$ of vertical
2-forms on $M$. We find the matrix of the endomorphism as the function system
$\{\,r_{\b\c\d\e}\,\}$ on the space of fibre bundle $B{\sf Q}$. Let
$\{\,R_{\b b\c c\d d\e h}\,\}$ be components of tensor $R$ on the space of
$\cal G$-structure, $\Omega$ be a vertical 2-form on $M$. Then
$R(\Omega)_{\b b\c c}=R_{\b b\c c\d d\e h}\,\Omega^{\d d\e h}$. But by (41),
$\Omega^{\d d\e h}=\om^{\d\e}\,G^{dh}$, $R(\Omega)_{\b b\c c}=r(\om)_{\b\c}\,G_{bc}$
and thus, $r(\om)_{\b\c}\,G_{bc}=R_{\b b\c c\d d\e h}\,G^{dh}\,\om^{\d\e}$,
hence, $r(\om)_{\b\c}=\frac1nG^{bc}\,G^{dh}\,R_{\b b\c c\d d\e h}\,\om^{\d\e}$. Thus,
\begin{equation}
r_{\b\c\d\e}=\frac1nG^{bc}\,G^{dh}\,R_{\b b\c c\d d\e h}.
\end{equation}
In view of symmetry properties of tensor $R$ and (56) we get
\procl{Proposition 15} Tensor $r$ has the following symmetry properties:
$$1)\,r_{\b\c\d\e}=-r_{\c\b\d\e};\quad 2)\,r_{\b\c\d\e}=-r_{\b\c\e\d};\quad
3)\,r_{\b\c\d\e}=r_{\d\e\b\c}.$$
\eprocl
\demo{Proof} Taking into account equality (56), we have:
\begin{eqnarray*}
1)\,r_{\b\c\d\e}\eq\frac1nR_{\b b\c c\d d\e h}\,G^{bc}\,G^{dh}
=-\frac1nR_{\c c\b b\d d\e h}\,G^{bc}\,G^{dh}=\\
\eq-\frac1nR_{\c c\b b\d d\e h}\,G^{cb}\,G^{dh}=-r_{\c\b\d\e}.
\end{eqnarray*}
Similarly for 2). Finally,
\begin{eqnarray*}
3)\,r_{\b\c\d\e}\eq\frac1nR_{\b b\c c\d d\e h}\,G^{bc}\,G^{dh}
=\frac1nR_{\d d\e h\b b\c c}\,G^{bc}\,G^{dh}=\\
\eq\frac1nR_{\d b\e c\b d\c h}\,G^{dh}\,G^{bc}=\frac1nR_{\d b\e c\b d\c h}\,
G^{bc}\,G^{dh}=r_{\d\e\b\c}.\qeds
\end{eqnarray*}
\edemo
Summarizing the above we get the following result:
\procl{Theorem 31} Riemann-Christoffel endomorphism of an $\hq$-manifold of
verti\-cal type induces a self-conjugate endomorphism of vertical 2-form module
of this manifold.\qed
\eprocl
\definition{Definition 25} Endomorphism $r:(\L^2)_V(M)\rightarrow(\L^2)_V(M)$
will be called {\it a twistor curvature\/} (or $t${\it -curvature}) {\it tensor\/}
of an $\hq$-manifold $M$.
\edefinition
The sense of the definition will be made clear in the following chapter. Note
that in the case $\dim M=4$ we have $r=R$, in particular, tensor $r$ components
in the frame $(p,e,J_1(e),J_2(e),J_3(e));\quad e\in T_p(M)$, $|e\|=1$,
coincides with tensor $R$ components in this frame.
The question naturally arises: When does endomorphism $r$ (in natural reprsentation)
preserve the module of self-dual (and anti-self-dual) forms of the manifold ?
We know that this module in natural representation corresponds to the Kaehler
module of the manifold, and thus, the question can be put in the following way:
When does Riemann-Christoffel endomorphism preserve the Kaehler module of an
$\hq$-manifold of vertical type ?
\definition{Definition 26} The structural bundle $\sf Q$ of an $\hq$-manifold
$M$ of vertical type will be called {\it Einsteinian\/} if on the space of
fibre bundle $B{\sf Q}$ over this manifold
\begin{equation}
g^{\b\d}r_{\b\c\d\e}=cg_{\c\e};\qquad c\in C^\infty\,(M).
\end{equation}
\edefinition
It is clear that if $\dim M=4$, the property of the structural bundle being
Einsteinian is equivalent to manifold $M$ being Einsteinian.
\procl{Theorem 32} Riemann-Christoffel endomorphism of an $\hq$-manifold of
verti\-cal type preserves the Kaehler module of the manifold iff the structural
bundle of the manifold is Einsteinian.
\eprocl
\demo{Proof} The assertion of the theorem is, evidently, equivalent to the
following: endomorphism $r$ preserves in natural representation module $\L^+(M)$
of self-dual forms on $\hq$-manifold $M$ of vertical type iff on $B{\sf Q}$
identities (57) are valid.  Let $\om\in\L^+(M)$. On the space $B({\sf Q})$
this 2-form is characterized by the function matrix
\begin{equation}
(\om_{\b\c})=\left(\begin{array}{rrrr}
0&-x&-y& -z\\x&0&-z&-\a y\\y&z&0&\a x\\z&\a y&-\a x&0\end{array}\right).
\end{equation}
Thus, $r(\L^+(M))\subset\L^+(M)\iff(\forall\,\om\in\L^+(M)\implies r(\om)\in\L^+(M))$, i.e.
\begin{equation}
\begin{array}{ll}
&1)\ (\a r_{01\b\c}+r_{23\b\c})\om^{\b\c}=0;\nonumber\\
&2)\ (\a r_{02\b\c}-r_{13\b\c})\om^{\b\c}=0;\\
&3)\ (r_{03\b\c}-r_{12\b\c})\om^{\b\c}=0.\nonumber
\end{array}
\end{equation}
Further, $\om^{\b\c}=g^{\b\d}\,g^{\c\z}\,\om_{\d\z}=\e(\b,\c)\om_{\b\c}$, and by (58),
$$(\om^{\b\c})=\left(\begin{array}{rrrr}
0&\a x&\a y&-z\\-\a x&0&-z&y\\-\a y&z&0&-x\\z&-y&x&0\end{array}\right).$$
Consequently, correlations (59) have the form
$$\begin{array}{rl}
&1)\ \ (\a r_{0101}+r_{2301})\a x+(\a r_{0102}+r_{2302})\a y-\\
&\quad-(\a r_{0103}+r_{2303})z-(\a r_{0112}+r_{2312})z+\\
&\quad+(\a r_{0113}+r_{2313})y-(\a r_{0113}+r_{2323})x=0.\end{array}$$
In view of the fact that these correlations must equivalently hold with respect
to $x,y,z$, we get:
\begin{equation}
\begin{array}{l}
r_{0101}-r_{2323}=0;\\r_{0102}+\a r_{2302}+\a r_{0113}+r_{2313}=0;\\
\a r_{0103}+r_{2303}+\a r_{0112}+r_{2312}=0.\end{array}\end{equation}
$$\begin{array}{rl}
&2)\ \ (\a r_{0201}-r_{1301})\a x+(\a r_{0202}-r_{2302})\a y-\\
&\quad-(\a r_{0203}-r_{1303})z-(\a r_{0212}-r_{1312})z+\\
&\quad+(\a r_{0213}-r_{1313})y-(\a r_{0223}-r_{1323})x=0;\end{array}$$
hence
\begin{equation}
\begin{array}{l}
r_{0201}-\a r_{1301}-\a r_{0223}+r_{1323}=0;\\
r_{0202}-r_{1313}=0;\\
\a r_{0203}-r_{1303}+\a r_{0212}-r_{1312}=0.\end{array}\end{equation}
$$\begin{array}{rl}
&3)\ \ (r_{0301}-r_{1201})\a x+(r_{0302}-r_{1302})\a y-\\
&\quad-(r_{0303}-r_{1203})z-(r_{0312}-r_{1212})z+\\
&\quad+(r_{0313}-r_{1213})y-(r_{0323}-r_{1223})x=0;\end{array}$$
hence
\begin{equation}
\begin{array}{l}
\a r_{0301}-\a r_{1201}-r_{0323}+r_{1223}=0;\\
\a r_{0302}-\a r_{1202}+r_{0313}-r_{1213}=0;\\
r_{0303}-r_{1212}=0.\end{array}\end{equation}
Comparing $(60_2)$ and $(61_1)$ we get:
\begin{equation}
r_{0102}+r_{1323}=0;\qquad r_{2302}+r_{0113}=0.
\end{equation}
Similarly, comparing $(60_3)$ and $(62_1), (61_3)$ and $(62_2)$, we get, respectively:
\begin{equation}
\begin{array}{ll}
\a r_{0112}+r_{0323}=0;\qquad&\a r_{0103}+r_{1223}=0;\\
\a r_{0212}-r_{0313}=0;\qquad&\a r_{0203}-r_{1213}=0.\end{array}
\end{equation}
Equalities (63) and (64) can be rewritted in the form
$$\begin{array}{lll}
r_{0102}+r_{3132}=0;\quad&r_{1013}+r_{2023}=0;\quad&-\a r_{1012}+r_{3032}=0;\\
r_{0103}-\a r_{2123}=0;\quad&-\a r_{2021}+r_{3031}=0;\quad&r_{0203}-\a r_{1213}=0;
\end{array}$$
or $g^{\b\d}\,r_{\b\c\d\e}=0\quad(\c\ne\e)$, i.e.
$g^{\b\d}\,r_{\b\c\d\e}=cg_{\c\e}\quad(\c\ne\e)$.
It is asserted that the system of the rest of the equalities
\begin{equation}
r_{0101}=r_{2323};\quad r_{0202}=r_{1313};\quad r_{0303}=r_{1212};
\end{equation}
is equivalent to equalities $g^{\b\c}\,r_{\b\c\d\c}=cg_{\c\c}$, i.e.
$$g^{\b\d}\,r_{\b 0\d 0}=-\a g^{\b\d}\,r_{\b 1\d 1}=-\a g^{\b\d}\,r_{\b 2\d 2}
=g^{\b\d}\,r_{\b 3\d 3},$$ i.e.
\begin{equation}
\begin{array}{l}
-\a r_{1010}-\a r_{2020}+r_{3030}=-\a r_{0101}+r_{2121}-\a r_{3131}\\
=-\a r_{0202}+r_{1212}-\a r_{3232}=r_{0303}-\a r_{1313}-\a r_{2323}.
\end{array}\end{equation}
Indeed, comparing the first and the forth expressions in (66), we get:
$$r_{1010}+r_{2020}=r_{1313}+r_{2323}.$$
On the other hand, comparing the second and the third expressions in (66), we get:
$$r_{1010}-r_{2020}=-r_{1313}+r_{2323},$$
hence,
$$r_{1010}=r_{2323},\quad r_{2020}=r_{1313}.$$
In view of the above and comparing the first and the second equalities in (66)
we find that $-\a r_{2020}+r_{3030}=r_{2121}-\a r_{2020}$, i.e. $r_{3030}=r_{1212}$.
Evidently, the inverse is also trouth. Combining these results we get that
$g^{\b\d}r_{\b\c\d\e}=cg_{\c\e}$. Inversely, from the above it follows that
validity of these conditions is equivalent to (63)-(65) that, as we have seen,
are  equivalent to (59) by (58), i.e. to module $\L^+(M)$ being invariant with
respect to endomorphism $r$.\qed
\edemo
\procl{Corollary} A 4-dimensional oriented pseudo-Riemannian manifold of
index 0 or 2 is Einsteinian iff its module of self-dual forms is invariant
with respect to Riemann-Chrictoffel endomorphism.\qed
\eprocl
The above Corollary shows that the proved Theorem 32 is a broad genera\-lization
of the known Atiyah-Hitchin-Singer theorem giving a criterion of 4-dimensional
Riemannian manifolds being Einsteinian in terms of self-dual forms [3] since
it generalizes this theorem in case of neutral pseudo-Riemannian metric itself.
\section{$t$-conformal-semiflat Manifolds}
The great importance of conformal-semiflat, i.e. self-dual or anti-self-dual,
manifolds in geometry and theoretical physics is explained by the fact that
the canonical almost complex structure on their twistor spaces is integrable
that allows to have twistor interpretation of Yang-Mills fields - instantons
and anti-instantons respectively, on the manifolds generalizing the classical
Ward theorem [3].
\par
In the previous chapter we introduced the notion of twistor curvature tensor
(endomorphism) of $\hq$-manifolds of vertical type. Now consider the trace-free
part of the endomorphism that we call  {\it a Weyl endomorphism.} Using it we
get a natural generalization of self-dual and anti-self-dual manifolds for the
case of $\hq$-manifolds of arbitrary dimension. We call them $t${\it-conformal-semiflat\/}
$\hq$-manifolds and show that their twistor curvature tensor is the algebraic
curvature tensor. The notion of twistor curvature of $\hq$-manifold is introduced
and $\hq$-manifolds of constant twistor curvature are studied. It is proved that
anti-self-dual $\hq$-manifolds, as well as generalized quaternionic-Kaehler
manifolds are $\hq$-manifolds of constant twistor curvature.
\subsection{Self-dual and anti-self-dual $\hq$-manifolds}
Let $M$ be an $\hq$-manifold of vertical type, $\sf Q$ be its structural
bundle, $r:(\L^2)_V(M)\rightarrow(\L^2)_V(M)$ be a twistor curvature tensor of
manifold $M$. Construct tensor $ric\in({\cal T}^0_2)_V(M)$ with components
$ric_{\b\e}=g^{\c\d}\,r_{\b\c\d\e}$ that we call {\it Ricci twistor tensor},
or {\it Ricci} $t${\it-curvature tensor.} By Theorem 32 a Riemann-Christoffel
endomorphism of an $\hq$-manifold of vertical type preserves the Kaehler module
of the manifold iff its Ricci twistor tensor is proportional to tensor $g$ of
fibre metric. The trace $\kappa=g^{\b\c}\,ric_{\b\c}$ of tensor $ric$ will be called
{\it a $t$-scalar curvature\/} of manifold $M$. By means of tensor $ric$ we
construct endomorphism $W$ of module $(\L^2)_V(M)$ with components
$$W_{\b\c\d\e}=r_{\b\c\d\e}+\frac12(ric_{\b\d}g_{\c\e}+ric_{\c\e}g_{\b\d}
-ric_{\b\e}g_{\c\d}-ric_{\c\d}g_{\b\e})+\frac\kappa6(g_{\c\d}g_{\b\e}-g_{\b\d}\,g_{\c\e})$$
generalizing the classical Weyl tensor of conformal curvature of a 4-dimensional
pseudo-Riemannian manifold. We call it {\it a Weyl endomorphism.} Evidently,
it has all symmetry properties of a twistor curvature tensor found in
Proposition 15. Besides, it is evident that $g^{\c\d}W_{\b\c\d\e}=0$. Similarly
to the proof of Theorem 32 for endomorphism $r$ we get for Weyl endomorphism:
\procl{Theorem 33} Submodules of self-dual and anti-self-dual forms
on an $\hq$-manifold of vertical type are invariant with respect to the action
of Weyl endomor\-phism.\qed
\eprocl
Thus, $W=W^++W^-$, where $W^\pm$ are endomorphism $W$ narrowings onto submodules
$\L^\pm(M)$, respectively.
\definition{Definition 27} An $\hq$-manifold $M$ of vertical type will be called
{\it self-dual\/} (resp., {\it anti-self-dual\/}) if $W^-=0$ (resp., $W^+=0)$.
An self-dual or anti-self-dual $\hq$-manifold will be called $t${\it-conformal-semiflat.}
\edefinition
>From the above we see that the notions generalize the corresponding notions of
4-dimensional Riemannian geometry that have become classical [5].
\par
Introduce  {\it self-duality index\/} $\xi$ of a $t$-conformal-semiflat
$\hq$-manifold equal to 1 for an self-dual and -1 for an anti-self-dual manifold.
\par
Let $M$ be an $\hq$-manifold, $\om\in\L^\pm(M)$, i.e. $*(\om)=\pm\om$. As we
have seen in the proof of Theorem 23, on the space of fibre bundle $B{\sf Q}$,
where $(g_{\b\c})=diag(1,-\a,-\a,1)$, these equalities will assume the form
$\om_{\hat\b\hat\c}=\pm\e(\b,\c)\,\om_{\b\c}$, respectively, where $\e(\b,\c)
= g_{\b\b}\,g_{\c\c};\quad\b,\c=0,1,2,3,4;\quad(\hat\b,\hat\c)$ is the pair
complementing the pair $(\b,\c)$ up to even permutation of indices (0,1,2,3).
It means that
\begin{equation}
(\om_{\b\c})=\left(\begin{array}{rrrr}
0&-x&-y&-z\\x&0&\mp z&\mp y\\y&\pm z&0&\pm\a x\\z&\pm\a y&\mp\a x&0
\end{array}\right);\
(\om^{\b\c})=\left(\begin{array}{rrrr}
0&\a x&\a y&-z\\-\a x&0&\mp z&\pm y\\-\a y&\pm z&0&\mp x\\z&\mp y&\pm x&0
\end{array}\right).\end{equation}
Now let $M$ be a $t$-conformal-semiflat $\hq$-manifold. If $M$ is self-dual,
then $W(\om)=0\quad(\om\in\L^-(M))$. If $M$ is anti-self-dual, then $W(\om)=0
\quad(\om\in L^+(M))$. By (67) these conditions will be written in the form
$(\a W_{\b\c01}+\xi W_{\b\c23})x+(\a W_{\b\c02}-\xi W_{\b\c13})y
+(\xi W_{\b\c12}-W_{\b\c03})z=0$. Since these correlations must equivalently
hold with respect to $x$, $y$ and $z$, we get that they are equivalent to the conditions
$$W_{\b\c01}=-\a\xi W_{\b\c23},\qquad W_{\b\c02}=\a\xi W_{\b\c13},\qquad
W_{\b\c03}=\xi W_{\b\c12}.$$
Write the correlations in view of symmetry properties of tensor $W$, that are
similar to symmetry properties of tensor $r$, found in Proposition 15:
$$\begin{array}{l}
W_{1002}=-W_{1332}=-\a\xi W_{0113}=\a\xi W_{0223};\\
W_{1003}=\a W_{1223}=-\xi W_{0112}=-\a\xi W_{0332};\\
W_{2003}=\a W_{2113}=\xi W_{0221}=\a\xi W_{0331};\\
W_{0101}=W_{2323}=-\a\xi W_{0123};\\W_{0202}=W_{1313}=-\a\xi W_{0231};\\
W_{0303}=W_{1212}=\xi W_{0312}.\end{array}$$
By definition of tensor $W$ the equalities will be rewritten, respectively, in the form:
\begin{equation}
\begin{array}{l}
1)\ r_{1002}+\a\xi r_{0113}=\frac12(ric_{12}-\xi\,ric_{03});\\
2)\ r_{1332}+\a\xi r_{0223}=\frac12(ric_{12}-\xi\,ric_{03});\\
3)\ r_{1003}+\xi r_{0112}=\frac12(ric_{13}-\a\xi\,ric_{02});\\
4)\ r_{1223}+\xi r_{0332}=\frac12(\xi\,ric_{02}-\a\,ric_{13});\\
5)\ r_{2003}-\xi r_{0221}=\frac12(ric_{23}+\a\xi\,ric_{01});\\
6)\ r_{2113}-\xi r_{0331}=\frac12(\xi\,ric_{01}-\a\,ric_{23});\\
7)\ r_{0101}+\a\xi r_{0123}=\frac12(\a\,ric_{00}-ric_{11})-\frac16\a\kappa;\\
8)\ r_{2323}+\a\xi r_{0123}=\frac12(\a\,ric_{33}-ric_{22})-\frac16\a\kappa;\\
9)\ r_{0202}+\a\xi r_{0231}=\frac12(\a\,ric_{00}-ric_{22})-\frac16\a\kappa;\\
10)\ r_{1313}+\a\xi r_{0231}=\frac12(\a\,ric_{33}-ric_{11})-\frac16\a\kappa;\\
11)\ r_{0303}-\xi r_{0312}=-\frac12(ric_{00}+ric_{33})+\frac16\kappa;\\
12)\ r_{1212}-\xi r_{0312}=\frac\a2(ric_{11}+ric_{22})+\frac16\kappa.
\end{array}\end{equation}
Besides, by definition of tensor $ric$,
\begin{equation}
\begin{array}{ll}
r_{1002}+r_{1332}=ric_{12};\qquad&r_{1003}-\a r_{1223}=ric_{13};\\
r_{2003}-\a r_{2113}=ric_{23};\qquad&-\a r_{0112}+r_{0332}=ric_{02};\\
r_{0113}+r_{0223}=-\a ric_{03};\qquad&-\a r_{0221}+r_{0331}=ric_{01}.
\end{array}\end{equation}
>From the above immediately follows
\procl{Theorem 34} A twistor curvature tensor of $t$-conformal-semiflat an
$\hq$-mani\-fold is an algebraic curvature tensor, tensor $W$ being its trace-free
part, i.e. the classical Weyl tensor with respect to tensor $r$.
\eprocl
\demo{Proof} By $(68_7)$, $(68_9)$ and $(68_{11})$ we have:
$r_{0123}+r_{0231}+r_{0312}=-\a\xi(r_{0101}+r_{0202})+\xi r_{0303}
+\frac12\xi(ric_{00} \a\,ric_{11}+ric_{00}-\a\,ric_{22}+ric_{00}+ric_{33})
-\frac12\xi\kappa=\xi(\a\,r_{0110}+\a\,r_{0220}-r_{0330})+\frac12\xi(ric_{00}
-\a\,ric_{11}-\a\,ric_{22}+ric_{33})+\xi\,ric_{00}\xi\kappa=-\xi\,ric_{00}
+\frac12\xi\kappa+\xi\,ric_{00}-\frac12\xi\kappa=0$.
By Proposition 15 we get that $r_{\b[\c\d\e]}=0$, and thus, $r$ is an algebraic
curvature tensor of module $(\L ^2)_V(M)$. For such tensor, as is well-known
[5], its trace-free part is uniquely defined by the formula that gives tensor $W$.\qed
\edemo
Consider some examples  of $t$-conformal-semiflat $\hq$-manifolds.
\definition{Definition 28} An $\hq$-structure that is at the same time a
$\pi AQ_\a$-structure, whose structural endomorphisms are parallel in
Riemannian con\-nection is called {\it a generalized hyper-Kaehler structure.}
\edefinition
\procl{Theorem 35} Any generalized hyper-Kaehler manifold of dimension
greater than 4 is $t$-conformal-semiflat.
\eprocl
\demo{Proof} As we will see below (Theorem 39), any generalized hyper-Kaehler
manifold $M$ is Ricci-flat. Besides, it, evidently, is a generalized
quaternionic-Kaehler manifold, and by (55) and Theorem 28, $R(\Omega)=0\quad
(\Omega\in(\L^2)_V(M))$. By definition of twistor curvature tensor $r=0$. Thus,
$ric=0$, and $W=0$. In particular, manifold $M$ is $t$-conformal-semiflat.\qed
\edemo
\procl{Corollary} All $\hq$-manifolds of dimension greater than 4, mentioned
in exam\-ple 2, section2.2, are $t$-conformal-semiflat.\qed
\eprocl
\subsection{Some properties of $t$-conformal-semiflat manifolds}
\procl{Theorem36} A twistor curvature endomorphism of an self-dual (resp.,
anti-self-dual) $\hq$-manifold $M$ transfers all anti-self-dual (resp., self-dual)
forms on $M$ into self-dual (resp., anti-self-dual) iff $M$ is a manifold  of
zero $t$-scalar curvature.
\eprocl
\demo{Proof} Endomorphism $r$ has the given property iff
$$1)\ \a r(\om)_{01}=\mp r(\om)_{23};\quad2)\ \a r(\om)_{02}=\pm r(\om)_{13};
\quad3)\ r(\om)_{03}=\pm r(\om)_{12};$$
$\om\in\L^\mp(M)$. Writing the equalities in view of (68) and (69) we get that
they are equivalent to the condition $ric_{00}-\a\,ric_{11}-\a\,ric_{22}+ric_{33}=0$,
or $\kappa=tr(ric)=g^{\b\c}\,ric_{\b\c}=0$.\qed
\edemo
\procl{Theorem 37} A twistor curvature endomorphism of an self-dual (resp.,
anti-self-dual)$\hq$-manifold $M$ transfers all vertical 2-forms on $M$ into
self-dual (resp., anti-self-dual) iff $M$ is manifold of zero Ricci $t$-curvature.
\eprocl
\demo{Proof} Let endomorphism $r$ transfers any vertical 2-form on $M$ into
self-dual (resp, anti-self-dual). In particular, any self-dual form is transferred
into self-dual (in the second case it is in view of endomorphism $r$ being
self-conjugated). By Theorem 32 $ric$ is a scalar endomorphism, and by Theorem 36
its trace is equal to zero. Thus, $ric=0$. Inversely, if $ric=0$, then by the
same Theorems, $r(\L^\pm(M))\subset\L^\pm(M),\quad r(\L^\mp(M))\subset\L^\pm(M)$,
and hence, $r(\L_V(M))\subset\L^\pm(M)$.\qed
\edemo
\procl{Corollary} Riemann-Christoffel endomorphism $R$ of a 4-dimensional
con\-formal-semiflat pseudo-Riemannian manifold of index 0 or 2 transfers all
anti-self-dual forms on the manifold into self-dual (or self-dual into anti-self-dual,
depending on the choice of orientation of manifold) iff the scalar curvature
of the manifold is zero. Here, $R$ transfers any 2-forms on manifold into
self-dual (or anti-self-dual) iff the manifold is Ricci-flat.\qed
\eprocl
Let $M$ be a self-dual (resp., anti-self-dual) $\hq$-manifold, $p\in M$,
$\om\in\L^\mp(M)$ is a non-isotropic eigenvector of its endomorphism $r$ of
twistor curvature in the point. Without loss of generality the vector norm can
be considered equal to either $\sqrt{-\a}$ or 1. In the former case we choose
orthonormalized basis $\{\id,J_1,J_2,J_3\,\}$ of the space ${\sf Q}_p(M)$ so
that $\om$ coincides with the pseudofundamental (resp., fundamental) form of
$\a$-quaternion $J_1$.
By (67), in this basis
$$(\om^{\b\c})=\left(\begin{array}{rrrr}
0&\a&0&0\\-\a&0&0&0\\0&0&0&\xi\\0&0&-\xi&0\end{array}\right).$$
Thus equality $r(\om)=\l\om$ will be rewritten in view of (68) and (69) in the form
$$\begin{array}{l}
1)\ \frac23(\a\,ric_{00}-ric_{11})-\frac13(\a\,ric_{33}-ric_{22})=\a\l;\\
2)\ ric_{12}+\xi\,ric_{03}=0;\\
3)\ ric_{13}+\a\xi\,ric_{02}=0;\\
4)\ -\frac13(\a\,ric_{00}-ric_{11})+\frac23(\a\,ric_{33}-ric_{22})=\a\l.
\end{array}$$
Note that the first and third of the equalities can be further rewritten in
the form, respectively:
$$\begin{array}{rcl}
2(\a\,ric_{00}-ric_{11})-(\a\,ric_{33}-ric_{22})\eq3\a\l;\\
-(\a\,ric_{00}-ric_{11})+2(\a\,ric_{33}-ric_{22})\eq3\a\l;\end{array}$$
hence,
$$\begin{array}{ll}
1)\ \a ric_{00}-ric_{11}=3\a\l;\qquad&2)\ \a ric_{33}-ric_{22}=3\a\l;\\
3)\ ric_{12}+\xi\,ric_{03}=0;\qquad&4)\ ric_{13}+\a\xi\,ric_{02}=0;\end{array}$$
or
$$\begin{array}{ll}
1)\ ric^0_0+ric^1_1=3\l;\qquad&2)\ ric^3_3+ric^2_2=3\l;\\
3)\ ric^2_1-\a\xi\,ric^3_0=0;\qquad&4)\ ric^3_1-\xi\,ric^2_0=0.\end{array}$$
In particular, $\kappa=\tr ric=6\l$.
Now, let $\om$ have a positive norm. Assume $\|\om\|=1$. Choose orthonormalized
basis $\{\id,J_1,J_2,J_3\,\}$ of the space ${\sf Q}_p(M)$ so that
$\om$ coincides with pseudofundamental (resp., fundamental) form of $\a$-quaternion
$J_3$. By (67), in the basis
$$(\om^{\b\c})=\left(\begin{array}{rrrr}
0&0&0&-1\\0&0&\xi&0\\0&-\xi&0&0\\1&0&0&0\end{array}\right).$$
In this case equality $r(\om)=\l\om$ by (68) and (69) will be rewritten in the form
$$\begin{array}{l}
1)\ \frac23(ric_{00}+ric_{33})+\frac13\a(ric_{11}+ric_{22})=\l;\\
2)\ \frac13(ric_{00}+ric_{33})+\frac23\a(ric_{11}+ric_{22})=-\l;\\
3)\ ric_{13}-\a\xi\,ric_{02}=0;\\
4)\ ric_{23}+\a\xi\,ric_{01}=0.\end{array}$$
The first two equalities can be rewritten in the form, respectively:
$$\begin{array}{rcl}
2(ric_{00}+ric_{33})+\a(ric_{11}+ric_{22})\eq3\l;\\
-(ric_{00}+ric_{33})-2\a(ric_{11}+ric_{22})\eq3\l;\end{array}$$
hence,
$$\begin{array}{ll}
1)\ ric_{00}+ric_{33}=3\l;\qquad&2)\ ric_{11}+ric_{22}=-3\a\l;\\
3)\ ric_{13}-\a\xi ric_{02}=0;\qquad&4)\ ric_{23}+\a\xi ric_{01}=0;
\end{array}$$ or
$$\begin{array}{ll}
1)\ ric^0_0+ric^3_3=3\l;\qquad&2)\ ric^1_1+ric^2_2=3\l;\\
3)\ ric^1_3+\xi\,ric^0_2=0;\qquad&4)\ ric^2_3-\xi\,ric^0_1=0.\end{array}$$
In particular, $\kappa=\tr(ric)=6\l$. Thus we have proved
\procl{Theorem 38} Eigenvalue of non-isotropic anti-self-dual (or self-dual)
form $\om$, which is eigenvector of twistor curvature endomorphism of an
self-dual (resp., anti-self-dual) $\hq$-manifold is equal to $\frac16\kappa$,
where $\kappa$ is a $t$-scalar curvature of the manifold.\qed
\eprocl
By Theorem 36 we get
\procl{Corollary 1} If an self-dual (resp., anti-self-dual) $\hq$-manifpld
admits at any its point a non-isotropic anti-self-dual (resp., self-dual) form
belowing to the kernel of the twistor curvature endomorphism, then the $t$-scalar
curvature of the man fold at this point is zero, and thus, the image of any
anti-self-dual (resp., self-dual) form at the given point is an self-dual
(resp., anti-self-dual) form in view of the endomorphism.\qed
\eprocl
By Theorem 37 we also get
\procl{Corollary 2} Let $M$ be an self-dual (resp., anti-self-dual)
$\hq$-manifold. If its twistor curvature endomorphism at any point of the
manifold vanishes at least one non-isotropic anti-self-dual (resp., self-dual)
form, them $M$ is a manifold of zero $t$-scalar curvature. If the endomorphism
turns to zero  any  anti-self-dual (resp., self-dual) form on $M$, then $M$
is a manifold of Ricci-zero $t$-curvature.\qed
\eprocl
\procl{Corollary 3} Let $M$ be a 4-dimensional self-dual (resp., anti-self-dual)
pseudo-Riemannian manifold of index 0 or 2. If its Riemann-Ghristoffel endomorphism
at every point of the  manifold vanishes at least one non-isotropic anti-self-dual
(resp., self-dual) form, then $M$ is a zero scalar curvature manifold. If the
endomorphism vanishes any anti-self-dual (resp., self-dual) form on $M$, then
$M$ is a Ricci-flat manifold.\qed
\eprocl
\subsection{Twistor curvature of $\hq$-manifold}
Let $M$ be an $\hq$-manifold of vertical type, $J\in\{\,{\sf T}\,\}$ be a
twistor on $M$, $\Omega(X,Y)=\la X,JY\ra$ be its Kaehler form.
\definition{Definition 29} Function
$$k(J)=\frac{R(\Omega,\Omega)}{\|\Omega\|^2}=\frac{\la R(\Omega),\Omega\ra}{\|\Omega\|^2}$$
is called {\it a twistor curvature\/} (or $t${\it-curvature}) {\it of manifold\/}
$M$ {\it in the direction of twistor J.}
\edefinition
Since by Theorem 26 $\Omega$ is a vertical 2-form on $M$, its components on the
space of $\sf G$-structure (by (41) and Theorem 23) will have the form
$\Omega_{\b b\c c}=\om_{\b\c}\,G_{bc}$, where $\{\,\om_{\b\c}\,\}$ are components
on the space of fibre bundle $B{\sf Q}$ of self-dual form $\om$ which is the
fundamental form of $\a$-quaternion $J$. Thus,
$$\begin{array}{rcl}
\la R(\Omega),\Omega\ra\eq-\frac12R_{\b b\c c\d d\z h}\Omega^{\d d\z h}\Omega^{\b b\c c}\\
\eq-\frac12R_{\b b\c c\d d\z h}\om^{\d\z}\,\om^{\b\c}G^{dh}G^{bc}
=\frac n2r_{\b\c\d\z}\om^{\b\c}\om^{\d\z}=n(r(\om),\om),\end{array}$$
and since $\|\Omega\|^2=n\om^2$, we have:
$$k(J)=\frac{(r(\om),\om)}{\|\om\|^2},$$
and thus, tensor $r$  completely defines a twistor curvature on $M$,
that explains the name of the tensor.
\definition{Definition 30} An $\hq$-manifold $M$ is called {\it a manifold of
pointwise constant twistor curvature\/} $c$ if
$$\forall\ J\in\{\,{\sf T}\,\}\implies\la R(\Omega),\Omega\ra
=c\|\Omega\|^2;\qquad c\in C^\infty\,(M).$$
If here $c=const$, $M$ is called {\it a manifold of globally constant
twistor curvature.}
\edefinition
\example{Example 1} Let $M^{4n}$ be an $\hq$-manifold of constant curvature $c$. Then
$$\begin{array}{rcl}
\la R(\Omega),\Omega\ra\eq-\frac12R_{\b b\c c\d d\z h}\,\Omega^{\b b\c c}\,\Omega^{\d d\z h}\\
\eq\frac c2(G_{\b b\d d}\,G_{\c c\z h}-G_{\c c\d d}\,G_{\b b\z h})\Omega^{\b b\c c}\,\Omega^{\d d\z h}\\
\eq\frac c2(g_{\b\d}\,g_{\c\z}\,G_{bd}\,G_{ch}-g_{\c\d}\,g_{\b\z}\,G_{cd}\,G_{bh})\om^{\b\c}\,G^{bc}\,\om^{\d\z}\,G^{dh}\\
\eq c\om_{\d\z}\om^{\d\z}\,G_{cd}\,G^{cd}=2c\|\Omega\|^2,\end{array}$$
i.e. $M$ is a manifold of globally constant twistor curvature.\eexample
\example{Example 2} Let $M^{4n}$ be a generalized quaternionic-Kaehler manifold,
$\Omega\in{\cal K}(M)$, $I$ be its corresponding twistor. Then at an arbitrary
point $p\in M$, by (47) and (52) that are, evidently true in any manifold dimension,
$$\begin{array}{rcl}
\la R(\Omega),\Omega\ra\eq\frac12\la R(\Omega)_{ij}\Omega^{ij}\ra
=\frac12\sum_{i,j=1}^{4n}\la R(\Omega)e_i,e_j\ra\Omega(e_i,e_j)\\
\eq-\frac{\a n}{2(n+2)}\sum_{i,j=1}^{4n}Ric(e_i,Ie_j)\la e_i,Ie_j\ra\\
\eq-\frac{\a n}{2(n+2)}\sum_{i,j=1}^{4n}Ric(e_i,e_j)\la e_i,e_j\ra\\
\eq-\frac{\a n}{2(n+2)}\sum_{i=1}^{4n}Ric(e_i,e_i)\|e_i\|^2\\
\eq-\frac{\a ns}{2(n+2)}=-\frac{\a s}{4(n+2)}\|\Omega\|^2,\end{array}$$
because
$$\begin{array}{rcl}
\|\Omega\|^{2}\eq\frac12\Omega_{ij}\,\Omega^{ij}=\frac12\sum_{i,j=1}^{4n}\la e_i,e_j\ra^2
=\frac12\sum_{i,j=1}^{4n}\la e_i,I^2e_j\ra^2\\
\eq\frac12\sum_{i,j=1}^{4n}\la e_i,e_j\ra^2=\frac12\sum_{i=1}^{4n}\|e_i\|^2
=2n\end{array}$$
\eexample
Here $\{\,e_1,\ldots,e_{4n}\,\}$ is orthonormalized basis adapted to twistor $I$.
Thus, any generalized quaternionic-Kaehler manifold is a manifold of constant
twistor curvature $c=-\frac{\a s}{4(n+2)}$, wnere $s$ is the scalar curvature
of the manifold. Moreover, if $\dim M>4$, by Theorem 29, $M$ is an Einsteinian
manifold, and thus, the manifold of globally constant twistor curvature
$c=-\frac{\a n\e}{n+2}$, where $\e$ is Einstein constant. In particular, all
quaternionic-Kaehler manifolds given in examples 2 and 3 of section 2.2, are
manifolds of constant twistor curvature.
\example{Example 3} Let $\{\,I,J\,\}$ is a generalized hyper-Kaehler structure,
i.e. $\pi AQ_\a$-structure generated by twistors $I$ and $J$ that are parallel
in Riemannian connection, on pseudo-Riemannian manifold $M$. Let $J_1=I$,
$J_2=J$, $J_3=I\circ J$. Then $\n X(J_\b)=0;\quad X\in{\sf X}(M)\quad\b=1,2,3$.
Thus, $[R(X,Y),{\cal J}]=0;\quad X,Y\in{\sf X}(M),\quad{\cal J}\in\{\,{\sf T}\,\}$.
In particular, $R({\cal I})\circ{\cal J}={\cal J}\circ R({\cal I});\quad
{\cal I},{\cal J}\in\{\,{\sf T}\,\}$,
i.e. $R({\cal I})$ belongs to centre of $\a$-quaternion algebra, hence,
$R({\cal I})=c\id$. In view of the fact that $R(\Omega)$ is skew-symmetric tensor,
where $\Omega$ is a Kaehler form of twistor $\cal I$, we get $R(\Omega)=0$,
and thus, $\la R(\Omega),\Omega\ra=0\quad(\Omega\in{\cal K}(M))$, i.e. $M$ is a
manifold of zero twistor curvature.
\eexample
Moreover, since the generalized hyper-Kaehler structure is, evidently, gene\-ralized
quaternionic-Kaehler, by the above example we get that if $\dim M>4$, then
$\e=0$, i.e. $M$ is a Ricci-flat manifold. If $\dim M=4$, then by above example
$s=0$, and since $R(\Omega)=0\quad(\Omega\in{\cal K}(M))$, $M$ is an Einsteinian
manifold (by the Corollary of Theorem 32), and thus, $Ric=0$, i.e. $M$ is
Ricci-flat, too. Thus, we have proved
\procl{Theorem 39} Any generalized hyper-Kaehler manifold is Ricci-flat.\qed
\eprocl
The above theorem generalizes the known Berger result of Ricci curvature of
hyper-Kaehler manifolds [6].
\par
Let $M$ be an arbitrary $\hq$-manifold of vertical type of pointwise constant
twistor curvature $c$. From the above it follows that this is equivalent to
$$\la r(\om),\om\ra=c\la\om,\om\ra;\qquad\om\in\L^+(M).$$
Polarizing the equality and in view of the endomorphism $r$ being self-conjugated,
we find that $\la r(\om),\v\ra=c\la\om,\v\ra;\quad\om,\v\in\L^+(M)$. Thus,
$r(\om)=c\om+\varphi;\quad\varphi\in\L^-(M)$. Evidently, the inverse is also
true. Thus, we get
\procl{Theorem 40} An $\hq$-manifold $M$ of vertical type is a manifold of
pointwise constant twistor curvature $c$ iff $(r-c\id)\om\in\L^-(M);\om\in
\L^+(M)$.\qed
\eprocl
Let $M$ be an $\hq$-manifold, $\om\in(\L^2)_V(M)$, $\{\,\om_{\b\c}\,\}$ are
components of tensor $\om$ on the space $B{\sf Q}$. We have: $\om=\om^++\om^-;
\quad\om^\pm\in\L^\pm(M)$. By (67) we get:
$$(\om^+_{\b\c})=\left(\begin{array}{rrrr}
0&x&y&z\\-x&0&z&\a y\\-y&-z&0&-\a x\\-z&-\a y&\a x&0\end{array}\right);\quad
(\om^-_{\b\c})=\left(\begin{array}{rrrr}
0&a&b&c\\-a&0&-c&-\a b\\-b&c&0&\a a\\-c&\a b&-\a a&0\end{array}\right);$$
Hence, $\om_{01}=x+a$; $\om_{02}=y+b$; $\om_{03}=z+c$; $\om_{23}=-\a x+\a a$;
$\om_{13}=\a y-\a b$; $\om_{12}=z-c$. Thus,
\begin{equation}
x=\frac12(\om_{01}-\a\om_{23});\qquad y=\frac12(\om_{02}+\a\om_{13});\qquad
z=\frac12(\om_{03}+\om_{12}).
\end{equation}
These correlations define projection $\pi:(\L^2)_V(M)\rightarrow\L^+(M)$ along $\L^-(M)$.
\par
Now let $M$ be an anti-self-dual $\hq$-manifold. Compute $\pi\circ r(\om);\quad
\om\in\L^+(M)$. We have: $r(\om)_{\b\c}=-r_{\b\c\d\z}\,\om^{\d\z}$. By (68) and
(69) it is easy to compute that
$$r(\om)_{01}-\a r(\om)_{23}=\frac13\kappa x;\quad r(\om)_{02}+\a r(\om)_{13}
=\frac13\kappa y;\quad r(\om)_{03}+r(\om)_{12}=\frac13\kappa z.$$
In view of (70), this means that
$$x(\pi\circ r(\om))=\frac16\kappa x(\om);\qquad y(\pi\circ r(\om))=\frac16\kappa y(\om);
\qquad z(\pi\circ r(\om))=\frac16\kappa z(\om).$$
(Here $x$, $y$ and $z$ are coordinate functions). Thus, $\pi\circ r(\om)=\frac16\kappa\om$
and $\pi\circ(r-\frac16\kappa\id)(\om)=0$, i.e. $(r-\frac16\kappa\id)\om\in\L^-(M)$.
By Theorem 40 we get:
\procl{Theorem 41} An anti-self-dual $\hq$-manifold is manifold of pointwise
constant twistor curvature $c=\frac16\kappa$, where $\kappa$ is a $t$-scalar
curvature of the manifold.\qed
\eprocl
\procl{Corollary 1} A 4-dimensional anti-self-dual pseudo-Riemannian manifold
of index 0 or 2 is a manifold of pointwise constant twistor curvature
$c=\frac16 s$, where $s$ is a scalar curvature of the manifold.\qed
\eprocl
\procl{Corollary 2} A 4-dimensional anti-self-dual pseudo-Riemannian manifold
whose self-dual form module is invariant with  respect to the Riemann-Christoffel
endomor\-phism is a manifold of globally constant twistor curvature.
\eprocl
\demo{Proof} This immediately follows from the manifold being Einsteinian (by
Theorem 32), in particular, the manifold of constant scalar curvature.\qed
\edemo
\section{Self-dual Geometry of 4-dimensional Kaehler Manifolds}
The relation of geometry of conformal-semiflat manifolds to geometry of
Ein\-steinian manifolds [5], as well as to twistor geometry [7] is of great
interest to the researchers. For example, the well-known Penrose-Atiyah-Hitchin-Singer
theorem [3] asserts that the canonical almost complex structure of the twistor
space of a 4-dimensional oriented Riemannian manifold $(M,g)$ is integrable iff
$(M,g)$ is conformal-semiflat. Hitchin showed that if $(M,g)$ is also a compact
Einsteinian manifold of positive scalar curvature, then it is isometric to $S^4$
or ${\bf C}P^2$ with standard metrics [8]. This result is also closely connected
with twistor geometry: Hitchin proved that a 4-dimensional oriented compact
Riemannian  manifold $(M,g)$ has a Kaehler twistor space iff $(M,g)$ is
conformally equivalent to $S^4$ or ${\bf C}P^2$ with their standard conformal
structures [9]. B.-Y.Chen [10] and (independently and by another method)
J.P.Bourguignon [11] and A.Derdzinski [12] got a classification of self-dual
compact Kaehler manifolds, and M.Itoh [13] got a classification of self-dual
Kaehler-Einsteinian manifolds. Besides, M.Itoh gave a complete characteristics
of compact anti-self-dual Kaehler manifolds [13].
\par
In the present chapter we received a complete classification of self-dual
genera\-lised Kaehler manifolds (both of classical and non-exceptional Kaehler
manifolods of hyperbolic type) of constant scalar curvature. It is proved that
a generalized Kaehler manifold is anti-self-dual iff its scalar curvature is
equal to zero. The results essentially generalized the above mentioned results
of N.Hitchin, J.P.Bourguignon, A.Derdzinski, B.-Y.Chen and M.Itoh.
\subsection{Generalized almost Hermitian structures}
\definition{Definition 31 [14]} We call {\it a generalized\/} (in the narrow
sense) {\it almost Hermitian\/} (or, $\ah-${it structure\/}) on a manifold $M$
a pair $\{\,g=\la\cdot,\cdot\ra,J\,\}$ of tensor field on $M$, where $g$ is a
pseudo-Riemannian metric, $J$ is an endomorphism of module ${\sf X}(M)$, such
that $J^2=\a\id;\quad\a=\pm 1$, and called {\it a structural operator}, or
{\it structural endomorphism.} Here, $\la JX,JY\ra=-\a\la X,Y\ra;\quad
X,Y\in{\sf X}(M)$. If $\a=-1$, the $\ah$-structure is called {\it an almost
Hermitian structure of classical type\/}, if $\a=1$, it is called {\it an
almost Hermitian structure of hyperbolic type\/} or, in other terms, {\it an
almost para-Hermitian srtucture}. In case $\a=1$, as is well-known, metric $g$
is neutral, i.e. it has a zero signature.
\edefinition
A 2-form $\Omega(X,Y)=\la X,JY\ra;\quad X,Y\in{\sf X}(M)$ called {\it a fundamental
structure form\/} is associated to every $\ah$-structure. The form is non-degenerate
and, thus, its exterior square generates the orientation of manifold $M$.
\definition{Definition 32} We call {\it a canonical orientation\/} on
$\ah$-manifold $(M,g,J)$ an orientation coherent to the form $-\a\Omega\land\Omega$.
\edefinition
We know [14] that giving an $\ah$-structure $(M,J,g)$ on $M$ equivalent to
giving $G$-structure $\sf G$ on $M$ with structural group $G=U(2,{\bf K}_\a)$,
where ${\bf K}_\a$ is the field $\bf C$ of complex numbers in the classical
case or the ring $\bf D$ of double numbers in the hyperbolic case. The space
elements of the $G$-structure are called {\it frames adapted to the $\ah$-structure},
or $A${\it-frames} [14]. Recall the construction of $A$-frames. Let
$(p,e_0,e_1,Je_0,Je_1)$ be an orthonormalized frame on manifold $M$ adapted to
structural endomorphism $J$, $\|e_0\|^2=1$, $\|e_1\|^2=-\a$. Then its corresponding
$A$-frame $(p,\e_0,\e_1,e_{\hat 0},e_{\hat 1})$ is defined as a frame in the
space $T_p(M)\otimes{\bf K}_\a$ with vectors $\e_a=\frac1{\sqrt 2}(e_a+iJ^3e_a),
\e_{\hat a}=\frac1{\sqrt 2}(e_a-iJ^3e_a);\quad a,b,c\ldots=0,1$; $i$ is the
imaginary unit of the ring ${\bf K}_a$. Note that $\o\e_a=\e_{\hat a}$,
$J\e_a=i\e_a$, $J\e_{\hat a}=-i\e_a$, where $X\rightarrow\o X$
is the operator of natural involution (of conjugation) in ${\bf K}_\a$-modul
$T_p(M)\otimes{\bf K}_\a$.
\par
Compute the matrix of metric tensor $g$ components on the space $\sf G$. In
the frame $(p,e_0,e_1,Je_0,Je_1)$ we have by definition: $(g_{ik})=diag(1,-\a,-\a,1)$.
Thus, in the $A$-frame $(p,\e_0,\e_1,\e_{\hat 0},\e_{\hat 1})$ we have:
$g_{0\hat 0}=\la\e_0,\e_{\hat 0}\ra=\frac12\la e_0+i\a Je_0, e_0-i\a Je_0\ra
=\frac12\{\la e_0,e_0\ra+\la e_0,e_0\ra\}=1=g_{\hat 00}$. Similarly,
$g_{1\hat 1}=g_{\hat 11}=-\a$. The rest of the metric tensor components are
equal to zero. Thus,
\begin{equation}
(g_{ij})=(g^{ij})=\left(\begin{array}{rrrr}
0&0&1&0\\0&0&0&-\a\\1&0&0&0\\0&-\a&0&0\end{array}\right)
\end{equation}
$(i,j=0,1,\hat0,\hat1)$. Evidently, $\det (g_{ij})=1$. On the other hand, compute
fundamental form of the structure $\Omega(X,Y)=\la X,JY\ra$:
$\Omega_{0\hat0}=\la\e_0,J\e_{\hat0}\ra=-i=-\Omega_{\hat00}$. Similarly,
$\Omega_{1\hat1}=-\Omega_{\hat11}=i\a$. The rest of the components of the form
$\Omega$ are equal zero. Thus, $\Omega=\Omega_{ij}\,\om^i\land\om^j
=-2i\om^0\land\om^{\hat0}+2i\a\om^1\land\om^{\hat1}$. Hence, in particular,
$$\Omega\land\Omega=-4\a\sqrt{\det g}\,\om^0\land\om^1\land\om^{\hat0}\land\om^{\hat1}$$
and, thus, $A$-frames are positively oriented with respect to the canonical
orientation in the classical case. In particular, $\eta_{01\hat0\hat1}=-\a$.
In view of this equality and by (71) we can specify the conditions of self-duality
and anti-self-duality of 2-forms on an $\ah$-manifold, given by correlations
(7), on the space of $G$-structure:
\procl{Lemma 4} Let $(M,g,J)$ be a 4-dimensional generalized almost Hermitian
mani\-fold, $\om\in\L^2(M)$. In this case
\begin{enumerate}
\item $\om\in\L^+(M)$ iff on the space $\sf G$
$$(\om_{jk})=\left(\begin{array}{rrrr}
0&x+iz&iy&0\\-x-iz&0&0&-\a y\\-iy&0&0&x-iy\\0&i\a y&-x+iz&0\end{array}\right);$$
\item $\om\in\L^-(M)$ iff on the space $\sf G$
$$(\om_{jk})=\left(\begin{array}{rrrr}
0&0&iy&x+iz\\0&0&-x+iz&i\a y\\-iy&x-iz&0&0\\-x-iz&-i\a y&0&0\end{array}\right).$$
\end{enumerate}
\eprocl
\demo{Proof} Let $\om\in\L^\pm(M)$. By (7), (71) and the given remarks about
the components of volume form $\eta_g$ we have:
$$\begin{array}{l}
1)\ \pm\om_{01}=\frac12\eta_{01ij}\,g^{ik}\,g^{jm}\,\om_{km}
=\eta_{01\hat0\hat1}\,g^{\hat00}\,g^{\hat11}\,\om_{01}=\om_{01};\\
2)\ \pm\om_{0\hat0}=\frac12\eta_{0\hat0ij}\,g^{ik}\,g^{jm}\,\om_{km}
=\eta_{0\hat01\hat1}\,g^{1\hat1}\,g^{\hat11}\,\om_{\hat11}=-\a\om_{1\hat1};\\
3)\ \pm\om_{0\hat1}=\frac12\eta_{0\hat1ij}\,g^{ik}\,g^{jm}\,\om_{km}
=\eta_{0\hat1\hat01}\,g^{\hat00}\,g^{1\hat1}\,\om_{0\hat1}=-\om_{0\hat1};
\end{array}$$
Similarly, $\pm\om_{1\hat0}=-\om_{1\hat0};\quad\pm\om_{\hat0\hat1}=\om_{\hat0\hat1}$.
Besides, since $\om$ is a real tensor, $\om_{\hat i\hat j}=\o\om_{ij}\quad
(\hat{\hat a}=a)$. In particular, $\o\om_{a\hat a}=-\om_{a\hat a}$. All the correlations
are, evidently, equivalent to the assertion of Lemma 4.\qed
\edemo
\definition{Definition 33} An $\ah$-srtucture $(g,J)$ is called {\it a
generalized Kaehler\/} (or, $\k-$) {\it structure}, if structural endomorphism
$J$ is parallel in Riemannian connection $\nabla$ of the manifold.
\edefinition
In case $\a=-1$ this notion coincides with the classical notion of Kaehler
structure, in case $\a=1$ it coincides with the notion of Kaehler structure of
hyperbolic type, or in other terms, para-Kaehler structure [15]. Manifolds
carring a para-Kaehler structure are also known as {\it Rashevskii fibre bundles}.
The name is explained by the existence on such manifolds, first studied by
P.K.Rashevskii of two isotropic totally integrable totally geodesic distributions
- eigendistributions of the structural endomorphism complementary to each other [16].
\par
Let $(M,g,J)$ be a $\k$-manifold. We know [14] that the first group of Cartan
equations of such manifold on the space of $G$-structure $\sf G$ has the form:
\begin{equation}
1)\ d\om^a=\om^a_b\land\om^b;\qquad 2)\ d\om_a=-\om^b_a\land\om_b;
\end{equation}
where $\om_a=\e(a)\om^{\hat a}=\e(a)\o\om^a,\quad\o\om^a_b=-\e(a,b)\om^b_a;
\quad\e(a)=g_{aa},\quad\e(a,b)=g_{aa}\,g_{bb}$. By exterior differentiation of
(72) and in view of basis forms being linearly independent we get the second
group  of structural equations of $\k$-structure:
\begin{equation}
d\om^a_b=\om^a_c\land\om^c_b+A^{ad}_{bc}\,\om^c\land\om_d;
\end{equation}
where $\{\,A^{ad}_{bc}\,\}$ is the function system on space $\sf G$ symmetric
by subscripts or superscripts and satisfying the system of differential equations
\begin{equation}
dA^{ad}_{bc}+A^{ad}_{hc}\,\om^h_b+A^{ad}_{bh}\,\om^h_c-A^{hd}_{bc}\,\om^a_h
-A^{ah}_{bc}\,\om^d_h=A^{ad}_{bch}\,\om^h+A^{adh}_{bc}\,\om_h;
\end{equation}
where $\{\,A^{ad}_{bch},A^{adh}_{bc}\,\}$ is the function system on space $\sf G$
symmetric by any pair of subscrips or superscripts and serving as components
of covariant differential in the Riemannian connection of tensor $A$ defined,
in view of (74), by the function system $\{\,A^{ad}_{bc}\,\}$ and called
{\it the tensor of holomorphic sectional curvature\/} [17].
Comparing (72) and (73) with structural equations of the pseudo-Riemannian
structure we have
$$d\om^i=\om^i_j\land\om^j,\qquad d\om^i_j
=\om^i_k\land\om^k_j+\frac12R^i_{jkm}\om^k\land\om^m,$$
where $\{\,R^i_{jkm}\,\}$ are components of Riemann-Christoffel tensor $R$.
It is easy to compute the explicit components on space $\sf G$ (i.e. in $A$-frame):
\begin{equation}
R_{\hat abc\hat d}=R_{b\hat a\hat dc}=-R_{\hat ab\hat dc}
=-R_{b\hat ac\hat d}=\e(a,d)A^{ad}_{bc}.
\end{equation}
The rest of the tensor components are equal to zero. Now we can compute the
components of the trace-free part of tensor $R$, the Weyl tensor $W$ of
conformal curvature. We know that tensors $R$ and $W$ are connected by the
correlation $(\dim M=4)$:
\begin{equation}
\begin{array}{rcl}
W_{ijkm}\eq R_{ijkm}+\frac12(r_{ik}\,g_{jm}+r_{jm}\,g_{ik}-r_{im}\,g_{jk}-\\
\eq-r_{jk}\,g_{im})+\frac s6(g_{im}\,g_{jk}-g_{ik}\,g_{jm}),
\end{array}\end{equation}
where $r_{ik}=g^{jm}\,R_{ijmk}$ are Ricci tensor components, $s=g^{ik}\,r_{ik}$
is the scalar curvature of the manifold. In view of (71) and (75) it is easy
to compute that in $A$-frame $g_{ab}=g_{\hat a\hat b}=0,\quad
g_{a\hat b}=g_{\hat ba}=\e(b)\d^{b}_{a};\quad r_{ab}=r_{\hat a\hat b}=0;\quad
r_{a\hat b}=r_{\hat ba}=-\e(b)A^{cb}_{ca}$.
Thus, (76) in {\it A}-frame will have the form:
$W_{ab\hat c\hat d}=\frac12(r_{a\hat c}\,g_{b\hat d}+r_{b\hat d}\,g_{a\hat c}
-r_{a\hat d}\,g_{b\hat c}-r_{b\hat c}\,g_{a\hat d})
+\frac s6(g_{a\hat d}\,g_{b\hat c}-g_{a\hat c}\,g_{b\hat d}),$ or
\begin{equation}
W_{ab}{}^{cd}=\frac12(r^c_a\,\d^d_b+r^d_b\,\d^c_a-r^d_a\,\d^c_b-r^c_b\,\d^d_a)
+\frac s6(\d^d_a\,\d^c_b-\d^c_a\,\d^d_b).
\end{equation}
Similarly, $W_{a\hat bc\hat d}=R_{a\hat bc\hat d}-\frac12(r_{a\hat d}\,g_{\hat bc}
+r_{\hat bc}\,g_{a\hat d})+\frac s6g_{\hat bc}\,g_{a\hat d}$, or
\begin{equation}
W_a{}^b{}_c{}^d=-A^{bd}_{ac}-\frac12(r^d_a\,\d^b_c+r^b_c\,\d^d_a)
+\frac s6\d^b_c\,\d^d_a,
\end{equation}
as well as correlations defined by symmetry properties of tensor $W$ (similar
to symmetry properties of Riemann-Christoffel tensor).
\subsection{4-dimensional self-dual $\k$-manifolds}
Let $M$ be a 4-dimensional oriented pseudo-Riemannian manifold of index 0 or 2,
$W$ be its Weyl tensor of conformal curvature. Consider $W$ as an endomorphism
of module $\L^2(M)$. It is known [5] that submodules $\L^\pm(M)$ are invariant
with respect to the endomorphism, and thus, $W=W^++W^-$, where $W^\pm$ are
restrictions of endomorphism $W$ onto submodules $\L^\pm(M)$, respectively.
\definition{Definition 34 ([3],[10])} A 4-dimensional oriented pseudo-Riemannian
manifold of index 0 (or 2) is called {\it self-dual\/} (respectively,
{\it anti-self-dual\/}) if $W^-=0$ (respectively, $W^+=0$) for it. An self-dual
or anti-self-dual manifold is called {\it conformal-semiflat.}
\edefinition
\example{Example} Let $(M,g,J)$ be a 4-dimensional generalized complex space
form, i.e. a generalized Kaehler manifold of constant holomorphic sectional
curvature. It is easy to check that a $\k$-manifold is a generalized complex
space form iff
$$A^{ad}_{bc}=c\,\tilde\d^{ad}_{bc},$$
where $c=const,\quad\tilde\d^{ad}_{bc}=\d^a_b\,\d^d_c+\d^a_c\,\d^d_b$. In this
case $r^a_b=\e(a)r_{\hat ab}=-A^{ca}_{cb}=-c\,\tilde\d^{ca}_{cb}=-3c\,\d^a_b,\quad
s=2r_{a\hat b}\,g^{a\hat b}=2r^a_a=-12c$, and by (77) and (78) $W_{ab}{}^{cd}
=-c\,\d^{cd}_{ab}$, where $\d^{cd}_{ab}=\d^c_a\,\d^d_b-\d^c_b\,\d^d_a;\quad
W_a{}^b{}_c{}^d=\e(b,d)R_{a\hat bc\hat d}+3c\,\d^d_a\,\d^b_c-2c\,\d^b_c\,\d^d_a
=-A^{bd}_{ac}+c\d^b_c\,\d^d_a=-c\,\tilde\d^{bd}_{ac}+c\d^b_c\,\d^d_a
=-c\d^b_a\,\d^d_c$. Thus, if $\v\in\L^-(M)$, then by Lemma 4
$$W(\v)_{ab}=-W_{ab}{}^{cd}\,\v_{cd}=0;\qquad
W(\v)^a{}_b=-W^a{}_b{}^c{}_d\,\v^d{}_c=-c\,\d^a_b\,\v^c{}_c=0,$$
i.e. $W(\om)=0$. Thus, we get
\procl{Theorem 42} A 4-dimensional generalized complex space form is an
self-dual manifold.\qed
\eprocl
\eexample
Now let $M$ be an arbitrary 4-dimensional self-dual generalized Kaehler manifold,
$\v\in\L^-(M)$. Then $W(\v)=0$, and by (78) $W(\v)^b{}_a=B^{bd}_{ac}\,\v^c{}_d=0$,
where $$B^{bd}_{ac}=A^{bd}_{ac}+\frac12(r^b_c\,\d^d_a+r^d_a\,\d^b_c)-\frac16\d^b_c\,\d^d_a.$$
In view of Lemma 4 this equality will be written in the form:
$(B^{b1}_{a1}-B^{b0}_{a0})iy+(\a B^{b0}_{a1}+B^{b1}_{a0})x+(\a B^{b0}_{a1}-B^{b1}_{a0})iz=0$.
In view of $x,y,z\in C^\infty\,(M)$ being arbitrary, we have:
\begin{equation}
B^{b1}_{a1}-B^{b0}_{a0}=0;\qquad B^{b0}_{a1}=B^{b1}_{a0}=0;\qquad (a,b=0,1).
\end{equation}
On the other hand,
\begin{equation}
B^{bc}_{ac}=A^{bc}_{ac}+\frac12(r^b_a+r^b_a)-\frac s6\,\d^b_a=-\frac s6d^b_a.
\end{equation}
Comparing (79) and (80) we get $B^{b0}_{a0}=B^{b1}_{a1}=-\frac s{12}\,\d^b_a$.
Thus, $B^{bd}_{ac}=-\frac s{12}\,\d^b_a\,\d^d_c$. In other word,
\begin{equation}
A^{bd}_{ac}=-\frac12(r^b_c\,\d^d_a+r^d_a\,\d^b_c)+\frac s6\,\d^b_c\,\d^d_a
-\frac s{12}\,\d^b_a\,\d^d_c.
\end{equation}
Alternating the equqlity by indices $b$ and $d$ we get:
\begin{equation}
r^b_c\,\d^d_a+r^d_a\,\d^b_c-r^b_a\,\d^d_c-r^d_c\,\d^b_a
=\frac s2\,\d^{bd}_{ac}.
\end{equation}
It is immediately checked that the above is equivalent to the condition
$r^a_a=\frac 12s$ that holds automatically in view of Ricci tensor being real.
On the other hand, symmetrize (81) by indices $b$ and $d$:
\begin{equation}
A^{bd}_{ac}=-\frac14(r^b_c\,\d^d_a+r^d_a\,\d^b_c+r^b_a\,\d^d_c
+r^d_c\,\d^b_a)+\frac s{24}\,\tilde\d^{bd}_{ac}.
\end{equation}
In view of the above remark the equality is equivalent to (81).
Now, note that (83) is equivalent to the following:
\begin{equation}
A^{bd}_{ac}=t^{(b}_{(a}\,\d^{d)}_{c)};
\end{equation}
where $t$ is a tensor; the brackets mean symmetrization by the enclosed indices.
Indeed, let (84) be true, i.e. $4A^{bd}_{ac}=t^b_a\,\d^d_c+t^d_c\,\d^b_a
+t^b_c\,\d^d_a+t^d_a\,\d^b_c$. Contracting the equality by indices $d$ and $c$ we get:
$$-4r^b_a=4t^b_a+tr(t)\,\d^b_a,\mbox{ and thus, }\qquad
t^b_a=-r^b_a-\frac14tr(t)\,\d^b_a.$$
Contracting the equality, we get $tr(t)=-\frac13s$, and thus
\begin{equation}
t^b_a=-r^b_a+\frac s{12}\,\d^b_a.
\end{equation}
Inversely, substituting (85) into (84) we get (83).
\par
Let (84) and the equivalent (83) hold. Substituting (82) and (83) into (77)
and (78), respectively, we get:
$$W_{ab}{}^{cd}=\frac s6\,\d^{cd}_{ab};\qquad W_a{}^b{}_c{}^d
=\frac s{12}\,\d^b_a\,\d^d_c,$$
that, in view of Lemma 4, implies $W(\om)=0\quad(\om\in\L^-(M))$, i.e. $M$ is
an self-dual manifold. Thus, we have proved
\procl{Theorem 43} A 4-dimensional generalized Kaehler manifold is self-dual
iff\newline $A^{bd}_{ac}=t^{(b}_{(a}\d^{d)}_{c)},$ where $t$ is a tensor, and,
necessarily, $t^b_a=-r^b_a+\frac1{12}s\,\d^b_a.$\qed
\eprocl
\definition{Definition 35} A generalized Kaehler manifold is called {\it
non-exceptional\/} if its Ricci endomorphism has at least one non-isotropic
eigenvector at every point of the manifold.
\edefinition
\example{Example} Any Kaehler manifold of classical type (with definite metric)
is non-exceptional since, in view of its Ricci endomorphism $r$ being self-conjugate, it has at every  its  point  an
orthonormalized basis of tangent space at this point consisting of eigenvectors
of the endomorphism. Moreover, since on the space of $G$-structure $\sf G$
components $r_{ab}$ and $r_{\hat a\hat b}$ of Ricci endomorphism are equal to
zero,the endomorphism commutes with the structural endomorphism $J$, and thus,
the mentioned basis can be choose adapted to $J$, i.e. having the form
$(e_0,e_1,Je_0,Je_1)$. A-basis generated by this basis will be the eigenbasis
of endomorphism $r$, as well as $t$, and eigenvalues of the endomorphosms will
be real.
\par
Moreover, endomorphisms $r$ and $t$ of any non-exceptional $\k$-manifold have
similar properties. Indeed, in view of the fact that for a $\k$-manifold $M$
there holds the equality $r\circ J=J\circ r$, endomorphism $r$ preserves
eigensubmodules of endomorphism $J$, in particular, it admits narrowing on
eigensubmodules $D^{\pm i}_J$ of the endomorphism with eigenvalues $\pm i$,
respectively. Let $e$ be a non-isotropic eigenvector of endomorphism $r$ at a
point of $M$. Without loss of generality it can be considered unitary.
Then $Je$ is also a non-isotropic eigenvector of this endomorphism and it
corresponds to the same eigenvalue $\tilde\l$. Thus, $\e_0=\frac1{\sqrt2}(e+iJ^3e)$
is an eigenvector of endomorphism $r$ with the same (real) eigenvalue. Comlementing
the pair $\{\,\e_0,\e_{\hat0}\,\},\quad\e_{\hat0}=\frac1{\sqrt2}(e-iJ^3e)$ up
to $A$-basis $\{\,\e_0,\e_1,\e_{\hat0},\e_{\hat1}\,\}$ we get that in basis
$\{\,\e_0,\e_1\,\}$ of submodule $D^i_J$ at this point endomorphism $r$ has matrix
$$(r)=\left(\begin{array}{cc}\tilde\l&a\\0&\tilde\mu\end{array}\right);
\qquad a,\tilde\mu\in{\bf K}_\a.$$
On the other hand, from (73) it follows that $\o A{}^{ad}_{bc}=A^{bc}_{ad}$, and
thus $\o r{}^a_b=r^b_a$. Therefore, $a=0$, $\tilde\mu\in{\bf R}$, i.e. endomorphism
$r|D^i_J$ is diagonalized in basis $\{\,\e_0,\e_1\,\}$. Moreover, since in view
of endomorphism $r$ being real $\o r{}^a_b=r^{\hat b}_{\hat a}$, $r|D^{-i}_J$ is
diagonalized in basis $\{\,\e_{\hat 0},\e_{\hat 1}\,\}$ and given by the same
matrix. Thus, endomorphism $r$ (as well as $t$) is diagonalized both in
$A$-basis $\{\,\e_0,\e_1,\e_{\hat 0},\e_{\hat 1}\,\}$ and in the corresponding
orthonormalized basis $\{\,e_0,e_1,Je_0,Je_1\,\}$.
\eexample
Now, let $(M,g,J)$ be a non-exceptional a $\k$-manifold. By the above we can
consider subbundle ${\sf G}_1$ of $G$-srtucture $\sf G$, consisting of
$A$-frames, where
\begin{equation}
(t^b_a)=\left(\begin{array}{cc}\l&0\\0&\mu\end{array}\right);\qquad
\l=\tilde\l-\frac s{12},\ \mu=\tilde{\mu}-\frac s{12}.
\end{equation}
Note that by Theorem 43, $\l+\mu=-\frac13s$. By (74) components of tensor $t$
on the space $\sf G$ satisfy the differential equations
$$dt^a_b+t^a_c\,\om^c_b-t^c_b\,\om^a_c=t^a_{bc}\,\om^c+t^{ac}_b\,\om_c,$$
where $t^a_{bc}=A^{ha}_{hbc}-\frac23A^{hr}_{hrc}\,\d^a_b$,
$t^{ac}_b=\e(c)\,A^{hac}_{hb}-\frac23A^{hrc}_{hr}\,\d^a_b$. In particular,
on the space of ${\sf G}_1$-structure, where $t^a_b=\l_b\,\d^a_b$, $\l_0=\l$,
$\l_1=\mu$, we have:
\begin{equation}
\begin{array}{l}
1)\ d\l=\l_c\,\om^c+\l^c\,\om_c;\qquad\l_c=t^0_{0c},\ \l^c=t^{0c}_0;\\
2)\ d\mu=\mu_{c}\,\om^c+\mu^c\,\om_c;\qquad\mu_c=t^1_{1c},\ \mu^c=t^{1c}_1;\\
3)\ (\l-\mu)\om^0_1=F_c\,\om^c+F^c\,\om_c;\qquad F_c=t^0_{1c},\ F^c=t^{0c}_1;\\
4)\ (\mu-\l)\om^1_0=G_c\,\om^c+G^c\,\om_c;\qquad G_c=t^1_{0c},\ G^c=t^{1c}_0.
\end{array}\end{equation}
We write structural equations of a generalized Kaehler structure on the space ${\sf G}_1$:
\begin{equation}
\begin{array}{l}
1)\ d\om^0=\om^0_0\land\om^0+\om^0_1\land\om^1;\\
2)\ d\om^1=\om^1_0\land\om^0+\om^1_1\land\om^1;\\
3)\ d\om_0=-\om^0_0\land\om_0-\om^1_0\land\om_1;\\
4)\ d\om_1=-\om^0_1\land\om_0-\om^1_1\land\om_1;\\
5)\ d\om^0_0=\om^0_1\land\om^1_0+\l\om^0\land\om_0+\frac14(\l+\mu)\om^1\land\om_1;\\
6)\ d\om^1_0=\om^1_0\land\om^0_0+\om^1_1\land\om^1_0+\frac14(\l+\mu)\om^1\land\om_0;\\
7)\ d\om^0_1=\om^0_0\land\om^0_1+\om^0_1\land\om^1_1+\frac14(\l+\mu)\om^0\land\om_1;\\
8)\ d\om^1_1=\om^1_0\land\om^0_1+\frac14(\l+\mu)\om^0\land\om_0+\mu\om^1\land\om_1.
\end{array}\end{equation}
By exterior differentiation of $(88_5)$ and by (88) we have:
\begin{eqnarray*}
\{-\frac12(\l-\mu)\om^1_0-\l^1\,\om_0-\frac14(\l_0+\mu_0)\om^1\,\}\land\om^0\land\om_1+\\
+\{\,\frac12(\l-\mu)\om^0_1 -\l_1\,\om^0+\frac14(\l^0+\mu^0)\om_1\,\}\land\om^1\land\om_0=0.
\end{eqnarray*}
In view of (87) and linear independence of basic forms we get $G^0=2\l^1$;
$G_1=\frac12(\l_0+\mu_0)$; $F_0=2\l_1$; $F^1=\frac12(\l^0+\mu^0)$, and thus,
\begin{eqnarray}
(\mu-\l)\om^1_0\eq G_0\,\om^0+\frac12(\l_0+\mu_0)\om^1+2\l^1\om_0+G^1\om_1;
\nonumber\\
(\mu-\l)\om^0_1\eq-2\l_1\,\om^0-F_1\,\om^1-F^0\,\om_0-\frac12(\l^0+\mu^0)\om_1.
\end{eqnarray}
By exterior differentiation of $(88_{6})$ and by (88) we have:
$$\begin{array}{l}
(\mu-\l)\om^1_0\land\om^0\land\om_0-(\l-\mu)\om^1_0\land\om^1\land\om_1+\\
+\frac12(\l_0+\mu_0)\om^0\land\om^1\land\om_0+\frac12(\l^1-\mu^1)\om^1\land\om_0\land\om_1=0.
\end{array}$$
In view of (89) and linear independence of basic forms we get
$G_0=0$; $G^1=0$; $\mu^1=3\l^1$, and thus,
\begin{equation}
\begin{array}{rcl}
(\mu-\l)\om^1_0\eq\frac12(\l_0+\mu_0)\om^1+2\l^1\,\om_0;\\
(\mu-\l)\om^0_1\eq-2\l_1\,\om^0-\frac12(\l^0+\mu^0)\om_1.\end{array}
\end{equation}
Similarly, by exterior differentiation of $(88_8)$ and by (88) we have:
$$\begin{array}{l}
\frac12(\l-\mu)\om^0_1\land\om^1\land\om_0-\frac12(\l-\mu)\om^1_0\land\om^0
\land\om_1+\frac14(\l_1+\mu_1)\om^1\land\om^0\land\om_0+\\
+\frac14(\l^1+\mu^1)\om_1\land\om^0\land\om_0+\mu_0\,\om^0\land\om^1\land
\om_1+\mu^0\om_0\land\om^1\land\om_1=0;\end{array}$$
hence, by (90) $\l_0=3\mu_0$, $\l^0=3\mu ^0$, $\mu_1=3\l_1$; $\mu^1=3\l^1$.
Thus, $d\l=3\mu_0\,\om^0+\l_1\,\om^1+3\mu ^0\,\om_0+\l^1\,\om_1;\quad
d\mu=\mu_0\,\om^0+3\l_1\,\om^1+\mu^0\,\om_0+3\l^1\,\om_1.$ In particular,
$$ds=-3(d\l+d\mu)=-12(\mu_0\,\om^0+\l_1\,\om^1+\mu^0\,\om_0+\l^1\,\om_1).$$
Let $M$ be a manifold of constant scalar curvature. Then it follow that
$\mu_0=\mu ^0=\l_1=\l^1=0$, i.e. $\l=const$, $\mu=const$.
Here (90) assumes the form:
$$1)\ (\mu-\l)\om^1_0=0;\qquad 2)\ (\mu-\l)\om^0_1=0.$$
Consider the possible cases.
\par
I. $\l\ne\mu$. Then $\om^1_0=\om^0_1=0$. By $(88_6)$ in this case $\l=-\mu$,
and (88) assumes the form:
$$\begin{array}{ll}
d\om^0=\om^0_0\land\om^0;\qquad&d\om^1=\om^1_1\land\om^1;\\
d\om_0=-\om^0_0\land\om_0;\qquad&d\om_1=-\om^1_1\land\om_1;\\
d\om^0_0=\l\,\om^0\land\om_0;\qquad&d\om^1_1=-\l\,\om^1\land\om_1.
\end{array}$$
In this case $M$ is locally holomorphically isometric to the product of
2-dimensional manifold $S^2_\l$ of constant curvature $\l$ by 2-dimensional
manifold $S^2_{-\l}$ of constant curvature $(-\l)$. Note that such manifold is
conformally flat [5].
\par
II. $\l=\mu$. Then $t^a_b=\l\,\d^a_b$, $A^{ad}_{bc}=\frac14\l\,\tilde\d^{ad}_{bc}$,
i.e. $M$ is a generalized complex space form. In this case, by [17], $M$ is
locally holomorphically isometric to one of the following manifolds:
\begin{enumerate}
\item Complex plane ${\bf C}^2$;
\item Complex projective plane ${\bf C}P^2$;
\item Complex hyperbolic plane ${\bf C}H^2$;
\item Double Euclidean plane ${\bf R}^2\Join{\bf R}^2$;
\item The space of null-pairs ${\bf R}P^2\odot{\bf R}P^2
=GL(3,{\bf R})/GL(2,{\bf R})\times GL(1,{\bf R})$
of real projective plane;
\end{enumerate}
equipped by a canonical Kaehler structure of classical (in cases 1,2,3) or
hyperbolic (in cases 4 and 5) types.
We get the following result:
\procl{Theorem 44} An self-dual nonexeptional generalized Kaehler manifold
of constant scalar curvature is locally holomorphically isometric to one of the
following mani\-folds: 1)\ ${\bf C}^2$; 2)\ ${\bf C}P^2$; 3)\ ${\bf C}H^2$;
4)\ $S^2_\l\times S^2_{-\l}$; 5)\ ${\bf R}^2\Join{\bf R}^2$;
6)\ ${\bf R}P^2\odot{\bf R}P^2$; equipped by a canonical Kaehler
structure of classical or, respectively, hyperbolic type.\qed
\eprocl
\procl{Corollary} Any self-dual compact Kaehler manifold of classical type
is holomor\-phically isometrically covered with one of the following manifolds:
1)\ ${\bf C}^2$; 2)\ ${\bf C}P^2$; 3)\ ${\bf C}H^2$; 4)\ $S^2_\l\times S^2_{-\l}$;
equipped by a canonical Kaehler structure.
\eprocl
\demo{Proof} It was shown by A.Derdzinski [12], that any such manifold is
locally symmetric, and thus, is a manifold of constant scalar curvature. To
complete the proof we apply Theorem 44 and note that all manifolds in the above
Corollary are simple-connected.\qed
\edemo
\remark{Remark} Self-dual generalized Kaehler manifolds of non-constant scalar
curvature do not admit such finite classification. A.Derdzinski showed [18]
that on ${\bf C}^2$ there exists a 2-parametric family of self-dual Kaehler
metrics of non-constant scalar curvature.
\eremark
\subsection{Self-duality and Bochner tensor}
The geometry of Kaehler manifolds is a complex analog of Riemannian geo\-metry:
such important notions of Riemannian geometry as sectional curvature, space
forms and many others have its complex counterpart which is of quite non-trivial
sense in geometry of Kaehlerian and, more generally, almost Hermitian manifolds.
One of such notions is Weyl's conformal curvature tensor which is the basic
object of study of conformal geometry. In 1949 S.Bochner introduced the complex
analog of the tensor for Kaehlerian manifolds [19]. The tensor introduced by
S.Bochner and called after his name, posseses all symmetry properties of
Riemann-Chrystoffel tensor and has meaning for arbitrary almost Hermitian
manifolds [20]. However, it has quite a complicated structure, and in spite of
a considerable number of works devoted to its studying we have comparatively
little information about its geometry.
\par
In the present section the Bochner curvature tensor of generalized Kaehlerian
manifolds is introduced. The structure of the tensor is studied. The main
result of this section states that Bochner curvature tensor of four-dimensional
generalized Kaehler manifold vanishes iff the manifold is self-dual. Hence it
is proved that Bochner-flat generalized Kaehler manifolds (i.e. generalized
Kaehler-Bochner manifolds) are natural generalization of self-dual generalized
Kaehler manifolds.
\definition{Definition 36} Let $M$ be $2n$-dimensional generalized Kaehler
manifold. The tensor $B$ of type (3,1) on $M$ defined by the equality
$$\begin{array}{rcl}
B(X,Y)Z\eq R(X,Y)Z+\la Y,Z\ra L(X)-\la X,Z\ra L(Y)+\la L(Y),Z\ra X-\\
&{}&-\la L(X),Z\ra Y-\la J^3(Y),Z\ra L\circ J(X)+\la J^3(X),Z\ra L\circ J(Y)-\\
&{}&-\la L\circ J^3(Y),Z\ra J(X)+\la L\circ J^3(X),Z\ra J(Y)+\\
&{}&+2\la L\circ J^3(X),Y\ra J(Z)+2\la J^3(X),Y\ra L\circ J(Z);
\end{array}$$
where $L=\frac1{2(n+2)}\,r+\frac s{8(n+1)(n+2)}\id$, is called {\it Bochner tensor
of manifold M.}
\edefinition
Let as compute components of Bochner tensor in $A$-frame:
$$\begin{array}{rcl}
B_{\hat abc\hat d}\eq\la B(\e_c,\e_{\hat d}),\e_b,\e_{\hat a}\ra
=\la R(\e_c,\e_{\hat d}),\e_b,\e_{\hat a}\ra+\la\e_{\hat d},\e_b\ra\la L(\e_c),\e_{\hat a}\ra-\\
&{}&-\la\e_c,\e_b\ra\la L(\e_{\hat d}),\e_{\hat a}\ra+\la L(\e_{\hat d}),\e_b\ra\la\e_c,\e_{\hat a}\ra
-\la L(\e_c),\e_b\ra\la\e_{\hat d},\e_{\hat a}\ra-\\
&{}&-la J^3(\e _{\hat d}),\e_b\ra\la L\circ J(\e_c),\e_{\hat a}\ra
+\la J^3(\e_c),\e_b\ra\la L\circ J(\e_{\hat d}),\e_{\hat a}\ra-\\
&{}&-\la L\circ J^3(\e_{\hat d}),\e_b\ra\la J(\e_c),\e_{\hat a}\ra
+\la L\circ J^3(\e_c),\e_b\ra\la J(\e_{\hat d}),\e_{\hat a}\ra+\\
&{}&+2\la L\circ J^3(\e_c),\e_{\hat d}\ra\la J(\e_b),\e_{\hat a}\ra
+2\la J^3(\e_c),\e_{\hat d}\ra\la L\circ J(\e_b),\e_{\hat a}\ra\\
\eq R_{\hat abc\hat d}+L^h_c\,\d^a_h\,\d^d_b+L^h_b\,\d^d_h\,\d^a_c
+L^h_c\,\d^a_h\,\d^d_b+
\\&{}&L^h_b\,\d^d_h\,\d^a_c+2L^h_c\,\d^d_h\,\d^a_b+2L^h_b\,\d^a_h\,\d^d_c\\
\eq R_{\hat abc\hat d}+L^a_c\,\d^d_b+L^d_b\,\d^a_c+L^a_c\,\d^d_b+L^d_b\,\d^a_c
+2L^d_c\,\d^a_b+2L^a_b\,\d^d_c\\
\eq R_{\hat abc\hat d}+2(L^a_c\,\d^d_b+L^d_b\,\d^a_c+L^d_c\,\d^a_b
+L^a_b\,\d^d_c).\end{array}$$
Analogously,
$$\begin{array}{rcl}
B_{\hat a\hat bcd}\eq R_{\hat a\hat bcd}
+\la\e_d,\e_{\hat b}\ra\la L(\e_c),\e_{\hat a}\ra
-\la\e_c,\e_{\hat b}\ra\la L(\e_d),\e_{\hat a}\ra+\\
&{}&+\la L(\e_d),\e_{\hat b}\ra\la\e_c,\e_{\hat a}\ra
-\la L(\e_c),\e_{\hat b}\ra\la\e_d,\e_{\hat a}\ra-\\
&{}&-\la J^3(\e_d),\e_{\hat b}\ra\la L\circ J(\e_c),\e_{\hat a}\ra
+\la J^3(\e_c),\e_{\hat b}\ra\la L\circ J(\e_c),\e_{\hat a}\ra-\\
&{}&-\la L\circ J^3(\e_d),\e_{\hat b}\ra\la J(\e_c),\e_{\hat a}\ra
+\la L\circ J^3(\e_c),\e_{\hat b}\ra\la J(\e_d),\e_{\hat a}\ra+\\
&{}&+2\la L\circ J^3(\e_c),\e_d\ra\la\e_{\hat b}),\e_{\hat a}\ra
+2\la J^3(\e_c),\e_d\ra\la L\circ J(\e_{\hat b}),\e_{\hat a}\ra\\
\eq 0+L^h_c\,\d^a_h\,\d^b_d-L^h_d\,\d^a_h\,\d^b_c+L^h_d\,\d^b_h\,\d^a_c
-L^h_c\,\d^b_h\,\d^a_d-\d^b_d\,L^h_c\,\d^a_h+\\
&{}&+\d^b_c\,L^h_d\,\d^a_h-L^h_d\,\d^b_h\,\d^a_c+L^h_c\,\d^b_h\,\d^a_d\\
\eq L^a_c\,\d^b_d-L^a_d\,\d^b_c+L^b_d\,\d^a_c-L^b_c\,\d^a_d-L^a_c\,\d^b_d
+L^a_d\,\d^b_c-L^b_d\,\d^a_c+L^b_c\,\d^a_d=0.\end{array}$$
Therefore,
\begin{equation}
B_{\hat abc\hat d}= B_{b\hat a\hat dc}=-B_{\hat ab\hat dc}
=-B_{b\hat ac\hat d}=A^{ad}_{bc}+ 8L^{(a}_{(b}\,\d^{d)}_{c)}.
\end{equation}
All other components of the tensor are evidently equal to zero.
\definition{Definition 37} A generalized Kaehler manifold is called
{\it Bochner-flat}, or {\it genera\-lized Kaehler-Bochner manifold}, if its
Bochner tensor is equal to zero identically.
\edefinition
Let $M$ be Bochner-flat generalized Kaehler manifold. In view of (91) it is
eqivalent to
$$A^{ad}_{bc}=-8L^{(a}_{(b}\,\d^{d)}_{c)},\qquad\mbox{or}\qquad
A^{ad}_{bc}=t^{(a}_{(b}\,\delta^{d)}_{c)},$$
where $t^a_b=-8L^a_b$. In view of Theorem 43 we get the following result:
\procl{Theorem 45} Four-dimensional generalized Kaehler manifold is
self-dual iff it is Bochner-flat.\qed
\eprocl
\subsection{4-dimensional anti-self-dual $\k$-manifolds}
Let $M$ be an arbitrary four-dimensional anti-self-dual generalized Kaehler
manifold, $\v\in\L^+(M)$. Then $W(\v)=0$, and by (77) and (78) we have:
\par
1)\ $W(\v)^b{}_a= B^{bd}_{ac}\,\v^c_d=0$. In view of Lemma 4, the equality will
be written in the form $(B^{b0}_{a0}+B^{b1}_{a1})iy=0$, and in view of
arbitrary choice of $y\in C^\infty\,(M)$, $B^{bc}_{ac}=0$, i.e.
$A^{bc}_{ac}+r^b_a-\frac16s\,\d^b_a=0$, and thus, $\frac16s\,\d^b_a=0$, i.e.
$s=0$. Inversely, if $s=0$, then $B^{bc}_{ac}=0$, and thus, $W(\v)^b{}_a=0$.
\par
2)\ $W(\v)_{ab}=W_{ab}{}^{cd}\,\v_{cd}=0$. In view of Lemma 4, the equality
will be rewritten in the form $W_{ab}{}^{01}(x+iz)=0$, and in view of arbitrary
choice of $x,z\in C^\infty\,(M)$, $W_{ab}{}^{01}=0$, and thus, $W_{ab}{}^{cd}=0$.
By (77) the above is equivalent to
\begin{equation}
r^c_a\,\d^d_b+r^d_b\,\d^c_a-r^c_b\,\d^d_a-r^d_a\,\d^c_b=\frac13\,\d^{cd}_{ab}.
\end{equation}
Contracting it by indices $b$ and $d$, we get: $2r^c_a+\frac12s\,\d^c_a
-2r^c_a=\frac13s\,\d^c_a$, hence, $s=0$. Inversely, if $s=0$, then (92) will
be written in the form $r^c_a\,\d^d_b+r^d_b\,\d^c_a-r^c_b\,\d^d_a-r^d_a\,\d^c_b=0$,
that is equivalent to $r^0_0+r^1_1=0$, i.e. $s=0$, and thus, $W(\v)_{ab}=0$. Hence,
$$W(\v)=0\iff s=0,$$
and we get the following result:
\procl{Theorem 46} A 4-dimensional generalized Kaehler manifold is
anti-self-dual iff it is manifold of zero scalar curvature.\qed
\eprocl
\remark{Remark} If $M$ is a 4-dimensional compact regular spinor manifold [7]
carrying a Kaehler structure of classic type, this result can be essentially
strengthened. Namely, in this case signature $\tau(M)$ of manifold $M$ is
computed by the formula [5]:
$$\tau(M)=\frac1{12\pi^2}\int\limits_M(\|W^+\|^2-\|W^-\|^2)\eta_g.$$
If $M$ is anti-self-dual and not conformal-flat, then it follows that $\tau(M)<0$.
But then $\hat A$-kind of the manifold [5], in 4-dimensional case equal to
$\frac1{16}\tau(M)$, is also negative, and by the Theorem of A.Lichneriwicz
[21] $M$ does not admit metrics of positive scalar curvature. By the Theorem
of Kazdan-Warner [22], the initial metric on $M$, being by Theorem 46 a metric
of zero scalar curvature, is Ricci-flat. If $M$ is conformal-flat, then since
it is a regular Kaehler surface, its first Betti number is zero, i.e. the
universal covering space is compact and, by Theorem 46, it has zero scalar
curvature, but it is impossible in view of Theorem 44.
\par
Thus, in this case $M$ is Ricci-flat. Inversely, if $M$ is Ricci-flat, then by
Theorem 46 it is anti-self-dual. Thus, we get the following result:
\procl{Theorem 47} A 4-dimensional compact regular spinor manifold carrying
any Kaehler structure of classical type is anti-self-dual iff it is Ricci-flat.\qed
\eprocl
An example of such manifold is a $K_3$-{\it surface}, i.e. a compact regular
surface with zero first Chern class. It is known [23] that such surfaces form
a special type in Kodaira classification of complex surfaces and they are well
studied. It is known, in particular, that all $K_3$-surfaces are diffeomorphic
and their Kaehler structures induced by Calabi-Yau metric are Kaehler-Einsteinian
structures [5]. Moreover, by the known result of Hitchin [8] any compact
anti-self-dual Einsteinian manifold is either flat or covered by a $K_3$-surface,
having a Calabi-Yau metric.
\eremark


\begin{thebibliography}{XX}
\bibitem {1} M. Berger, {\it Remarques sur le Groupe d'Holonomie des Vari\'et\'es
Rieman\-niennes}, C. R. Acad. Sci. Paris {\bf 262} (1966), 316-318.
\bibitem {2} S.M. Salamon, {\it Quaternionic Kaehler Manifolds}, Invent. Math.
{\bf 67} (1982), 143-171.
\bibitem {3} M.F. Atiyah, N.J. Hitchin, I.M. Singer, {\it Self-duality in
four-dimensional Riemannian Geometry}, Proc. Roy. London {\bf A362} (1978),
425-461.
\bibitem {4} S. Sternberg, {\it Lectures on Differential Geometry} (1970),
Mir, Moscow (Russian).
\bibitem {5} A. Besse, {\it Einsteinian Manifolds}, {\bf I-II} (1990),
Mir, Moscow, (Russian).
\bibitem {6} M. Berger, {\it Sur le Groupes d'Holonomie des Vari\'et\'es
\`a Connexion Affine et des Vari\'et\'es Riemanniennes}, Bull. Soc. Math.
France {\bf 83} (1955), 279-330.
\bibitem {7} R. Penrose, {\it The Twistor Programme}, Math. Phys. Repts.
{\bf 12} (1977), 65-76.
\bibitem {8} N.J. Hitchin {\it On Compact Four-dimensional Einstein manifolds},
J. Diff. Geom. {\bf 9} (1974), 435-442.
\bibitem {9} N.J. Hitchin, {\it Kahlerian Twistor Spaces}, Proc. London
Math. Soc. {\bf 43} (1981), 133-150.
\bibitem {10} B.-Y. Chen, {\it Some Topological Obstructions to Bochner-Kaehler
Metrics and their Appli\-cations} J. Diff. Geom. {\bf 13} (1978), 547-558.
\bibitem {11} J.-P. Bourguignon, {\it Les Vari\'et\'es de Dimension 4 a Signature
non Nulle dont la Courbure est Harmonique sont d'Einstein}, Invent. Math.
{\bf 63} (1981), 263-286.
\bibitem {12} A. Derdzinski, {\it Self-dual Kahler Manifolds and Einstein
Manifolds of Dimensional Four}, Comp. Math. {\bf 19} (1983), 405-433.
\bibitem {13} M. Itoh, {\it Self-duality of K\"ahler Surfaces}, Comp. Math.
{\bf 51} (1984), 265-273.
\bibitem {14} V.F. Kirichenko, {\it Generalized Hermitian Geometry Methods in
Theory of Almost Contact Manifolds}, in book Itogi Nauki i Tehniki: Problemy Geometrii
{\bf 18} (1986), 25-72 (Russian).
\bibitem {15} P. Libermann, {\it Sur le Probl\`em\'e d'Equivalence de Certaines
Structure Infinitesimales}, Ann di Mathematica {\bf 36} (1951), 247-261.
\bibitem {16} A.P. Shirokov, {\it On One Space Class over Algebras}, Izv.
VUZov. Matem. {\bf 1} (1961), 163-170 (Russian).
\bibitem {17} V.F. Kirichenko, {\it Axioms of Holomorphic Planes in Generalized
Hermitian Geometry}, Dokl. Akad. Nauk SSSR {\bf 260} (1981), 795-799 (Russian).
\bibitem {18} A. Derdzinski, {\it Examples of Kaehler-Einstein Self-dual Metrics
on Comp\-lex Plane},in book Four-dimensional Riemannian Geometry. Arthur Besse
Seminar 1978-1979. Mir, Moscow (1985), 280-289 (Russian).
\bibitem {19} S. Bochner, {\it Curvature and Betti Numbers, II} Ann. Math.
{\bf 50} (1949), 77-93.
\bibitem {20} S. Tachibana, {\it On Bochner Curvature Tensor}, Not. Sci.
Ochanomidzu Univ. {\bf 18} (1967), 15-19.
\bibitem {21} A. Lichnerowicz, {\it Spineur Harmoniques}, C. R. Acad. Sci.
Paris {\bf 257} (1963), 7-9.
\bibitem {22} J.L. Kazdan, F.W. Warner, {\it A Direct Approach to the
Determination of Gaussian and Scalar Curvature Function}, Invent. Math.
{\bf 28} (1975), 227-230.
\bibitem {23} K. Kodaira, {\it On the Structure of Compact Complex Analytic
Surfaces, I} Amer. J. Math. {\bf 86} (1964), 751-798.
\end{thebibliography}
\end{document}